\documentclass[12pt]{article}
\usepackage{amsfonts}

\def\hybrid{\topmargin 0pt      \oddsidemargin 0pt
        \headheight 0pt \headsep 0pt
        \voffset=-0.5cm
        \textwidth 6.5in        
        \textheight 9in         
        \marginparwidth 0.0in
        \parskip 5pt plus 1pt   \jot = 1.5ex}
\catcode`\@=11
\def\marginnote#1{}

\newcount\hour
\newcount\minute
\newtoks\amorpm
\hour=\time\divide\hour by60
\minute=\time{\multiply\hour by60 \global\advance\minute by-\hour}
\edef\standardtime{{\ifnum\hour<12 \global\amorpm={am}%
        \else\global\amorpm={pm}\advance\hour by-12 \fi
        \ifnum\hour=0 \hour=12 \fi
        \number\hour:\ifnum\minute<10 0\fi\number\minute\the\amorpm}}
\edef\militarytime{\number\hour:\ifnum\minute<10 0\fi\number\minute}

\def\draftlabel#1{{\@bsphack\if@filesw {\let\thepage\relax
   \xdef\@gtempa{\write\@auxout{\string
      \newlabel{#1}{{\@currentlabel}{\thepage}}}}}\@gtempa
   \if@nobreak \ifvmode\nobreak\fi\fi\fi\@esphack}
        \gdef\@eqnlabel{#1}}
\def\@eqnlabel{}
\def\@vacuum{}
\def\draftmarginnote#1{\marginpar{\raggedright\scriptsize\tt#1}}
\def\draftlabel#1{{\@bsphack\if@filesw {\let\thepage\relax
   \xdef\@gtempa{\write\@auxout{\string
      \newlabel{#1}{{\@currentlabel}{\thepage}}}}}\@gtempa
   \if@nobreak \ifvmode\nobreak\fi\fi\fi\@esphack}
        \gdef\@eqnlabel{#1}}
\def\@eqnlabel{}
\def\@vacuum{}
\def\draftmarginnote#1{\marginpar{\raggedright\scriptsize\tt#1}}

\def\draft{\oddsidemargin -.5truein
        \def\@oddfoot{\sl preliminary draft \hfil
        \rm\thepage\hfil\sl\today\quad\militarytime}
        \let\@evenfoot\@oddfoot \overfullrule 3pt
        \let\label=\draftlabel
        \let\marginnote=\draftmarginnote
   \def\@eqnnum{(\theequation)\rlap{\kern\marginparsep\tt\@eqnlabel}%
\global\let\@eqnlabel\@vacuum}  }

\def\numberbysection{\@addtoreset{equation}{section}
        \def\theequation{\thesection.\arabic{equation}}}

\def\underline#1{\relax\ifmmode\@@underline#1\else
        $\@@underline{\hbox{#1}}$\relax\fi}

\def\titlepage{\@restonecolfalse\if@twocolumn\@restonecoltrue\onecolumn
     \else \newpage \fi \thispagestyle{empty}\c@page\z@
        \def\thefootnote{\fnsymbol{footnote}} }

\def\endtitlepage{\if@restonecol\twocolumn \else  \fi
        \def\thefootnote{\arabic{footnote}}
        \setcounter{footnote}{0}}  
\relax

\numberbysection
\hybrid

\def\beq{\begin{equation}}
\def\eeq{\end{equation}}
\def\bea{\begin{eqnarray}}
\def\eea{\end{eqnarray}}
\def\p{\partial}
\def\G{\Gamma}
\def\g{\gamma}
\def\s{\sigma}

\def\C{{\cal C}}

\def\a{\alpha}
\def\b{\beta}
\def\e{\varepsilon}
\def\l{\lambda}

\def\A{{\cal A}}

\def\D{{\cal D}}
\def\F{{\cal F}}

\def\L{{\cal L}}

\def\O{{\cal O}}
\def\P{{\cal P}}

\def\SP{{\cal S}}
\def\dim{{\rm dim}}
\def\res{{\rm res}}

\def\TC{{\cal T}}
\def\wt{\widetilde}
\def\wh{\widehat}

\def\k{\mu}
\def \matrix #1 {\left(\begin{array}{cc} #1 \end{array}\right)}
\def\mod{{\ {\rm mod\ }}}
\newtheorem{theo}{Theorem}[section]
\newtheorem{prop}[theo]{Proposition}
\newtheorem{cor}[theo]{Corollary}
\newtheorem{lem}[theo]{Lemma}
\newtheorem{rem}[theo]{Remark}

\def\square{\hfill
{\vrule height6pt width6pt depth1pt} \break \vspace{.01cm}}
\def\bpf{\hfill\break{\it Proof.} }
\def\epf{\square}

\def\bbZ{{\mathbb Z}}
\def\bbC{{\mathbb C}}
\def\bbP{{\mathbb P}}
\def\bbN{{\mathbb N}}

\begin{document}
\begin{titlepage}
\title{Integrable discrete Schr\"odinger equations and
a characterization of Prym varieties by a pair of quadrisecants}
\author{S.Grushevsky
\thanks{Department of Mathematics, Princeton University,
Princeton, NJ 08544, USA; e-mail:
sam@math.princeton.edu. Research is supported in part by National
Science Foundation under the grant DMS-05-55867.}
\and I.Krichever
\thanks{Columbia University, New York, USA and Landau Institute for
Theoretical Physics, Moscow, Russia; e-mail:
krichev@math.columbia.edu. Research is supported in part by National
Science Foundation under the grant DMS-04-05519.}}

\date{}

\maketitle

\begin{abstract}
We prove that Prym varieties are characterized geometrically by the
existence of a symmetric pair of quadrisecant planes of the
associated Kummer variety. We also show that Prym varieties are
characterized by certain (new) theta-functional equations.
For this purpose we construct and study a
difference-differential analog of the Novikov-Veselov hierarchy.
\end{abstract}

\end{titlepage}
\section{Introduction}

An involution $\s:\G\longmapsto \G$ of a smooth algebraic curve $\G$
induces an involution $\s^*: J(\G)\longmapsto J(\G)$ of the Jacobian
of the curve. The odd subspace under this involution, i.e. the set
of $z\in J(\G)$ such that $\s^*(z)=-z$, as the kernel of a
homomorphism of abelian varieties, is the sum of a lower-dimensional
abelian variety, called the Prym variety (the connected component of
zero in the odd subspace), and a finite group. The Prym variety
naturally has a polarization induced by the principal polarization
on the Jacobian. However, this polarization is not principal, and
the Prym variety admits a natural principal polarization if and only
if $\s$ has at most two fixed points on $\G$ --- this is the case we will concentrate on.

The problem of
characterizing the locus $\P_g$ of Pryms of dimension $g$ in the
moduli space $\A_g$ of all principally polarized abelian varieties
(ppav's) is well-known and has attracted a lot of interest over the
years. In some sense the Prym varieties may be geometrically the
easiest to understand ppavs beyond Jacobians, and one
could hope that studying them would be a first step towards
understanding the geometry of more general abelian varieties as
well.

Recently, the Prym varieties for the case of an involution with two
fixed points were characterized in \cite{kr-prym} by one of the
authors by the property of the theta function
satisfying a certain partial differential equation, coming from the
theory of integrable Schr\"odinger equations. The main goal of this
paper is to give a geometric (and equivalent theta-functional)
characterization of Prym varieties corresponding to involutions
with no fixed points (the Prym varieties of involutions with two
fixed points arise as a degeneration of this case).

Let $B$ be an indecomposable complex symmetric matrix with positive
definite imaginary part.
It defines an indecomposable ppav $X:=\bbC^g/\Lambda$, where
$\Lambda:=\bbZ^g+B\bbZ^g\subset \bbC^g$. The Riemann theta function is given by the
formula
$$
\theta(B,z):=\sum\limits_{m\in \bbZ^g} e^{2\pi i(z,m)+\pi
i(B m,m)},\ \ (z,m)=m_1z_1+\ldots+m_gz_g\,
$$
for $z\in\bbC^g$. The theta functions of the second order are defined by the formula
$$
 \Theta[\e](B,z):=\sum\limits_{m\in\bbZ^g} e^{2\pi i\left(2m+\e,z\right)+\pi i
 \left(2m+\e,B(m+\frac{\e}{2})\right)}
$$
for $\e\in (\bbZ/2\bbZ)^g$. The Kummer
variety $K(X)$ is then defined as the image of the Kummer map
$$
K:z\longmapsto \{\Theta[\e](z)\}_{{\rm all\ }\e
\in(\bbZ/2\bbZ)^g}\in \bbP^{2^g-1}.
$$
A projective $(m-2)$-dimensional plane
$\bbP^{m-2}\subset\bbP^{2^g-1}$ intersecting $K(X)$ in at least $m$ points
is called an $m$-secant of the Kummer variety.

The Kummer images of Jacobians of curves were shown to admit many
trisecant lines (see \cite{fay}). It was then shown by Gunning
\cite{gun1} that the existence of a one-dimensional family of
trisecants in fact suffices to characterize Jacobians among all
ppavs. Welters, inspired by the Gunning's theorem and the Novikov's
conjecture proved later by Shiota \cite{shiota}, formulated in
\cite{wel1} the following conjecture: {\it If $K(X)$ has a
trisecant, and $X$ is indecomposable, then $X$ is a Jacobian}, and
this was recently proved by the second-named author in
(\cite{kr-schot,kr-tri}).

Prym varieties possess generalizations of some properties that
Jacobians possess. In \cite{bd} Beauville and Debarre, and in
\cite{fay2} Fay showed that the Kummer images of Prym varieties
admit many quadrisecant planes. Similarly to the Jacobian case, it
was then shown by Debarre in \cite{deb} that the existence of a
one-dimensional family of quadrisecants characterizes Pryms.
However, Beauville and Debarre in \cite{bd} showed that the
existence of a single quadrisecant plane to the Kummer variety does
not characterize Pryms.

In this paper we prove that the Prym varieties are characterized by
the existence of a {\it symmetric pair} of quadrisecants of the
corresponding Kummer variety --- i.e. of two different 2-planes each
intersecting the Kummer variety in 4 points, such that the points of
secancy for the two planes are related in some precise way. We
deduce this from a characterization by some theta-functional
equations, and study the associated discrete Schr\"odinger equations
along the way.

The structure of this work is as follows. In section 2 we prove that (the
Kummer images of) Prym varieties have symmetric pairs of quadrisecants.
This is done via an algebro-geometric construction of difference potential
Schr\"odinger operators that play a crucial role in all our further
considerations. Our construction is a discrete analog of the well-known
Novikov-Veselov algebro-geometric construction from \cite{nv} of a
potential for the two-dimensional Schr\"odinger operators. The latter is a
reduction of a more general construction of {\it finite-gap on one energy
level Schr\"odinger operators in a magnetic field} first proposed in
\cite{dkn} and based on the concept of the Baker-Akhiezer functions
introduced in \cite{kr1,kr2}.

The Baker-Akhiezer functions are uniquely defined by giving the data
of an algebraic curve with fixed local coordinates in the
neighborhoods of marked points, and by a divisor of its poles away
from the marked points. These functions are not algebraic --- they
have essential singularities at marked point. To the authors' surprise in
the case
of unramified double covers the very same arguments that led Novikov
and Veselov to a proof that certain constraints on the
algebro-geometric spectral data of \cite{dkn} are sufficient for
potential reduction show that the poles divisor of the
Baker-Akhiezer functions associated to unramified double covers of
algebraic curves is of degree {\it less} than one would expect from
the general theory. It turned out that this unexpected observation
is equivalent to some well-known results of Mumford on Prym
varieties \cite{mum1}.

In section 3 we introduce a discrete analog of the Novikov-Veselov
hierarchy and study its properties. It is a set of
difference-differential equations describing integrable
deformations of potential difference Schr\"odinger operators. In
section 4 of the paper, following the lines of works
\cite{kr-schot,kr-prym,kr-tri} we construct a wave solution, and then in
section 5 we finish the proof of our main result on the
characterization of the Prym varieties of unramified covers. It is
necessary to emphasize that, unlike the Jacobian case, the Prym
variety remains compact under certain degenerations of the curve.
No characterization of Prym varieties given in terms of the period
matrix of the Prym differentials can single out the possibility of
such degenerations.

\begin{theo}[Main theorem]
An indecomposable principally polarized abelian variety $(X,\theta)\in\A_g$
is (in the closure of) the locus $\P_g$ of Prym varieties of
unramified double covers if and only if there exist vectors
$A,U,V, W\in\bbC^g$ representing distinct points in $X$, none of them points of order two, and constants
$c_1,c_2,c_3,w_1,w_2,w_3\in\bbC$ such that one of the following
equivalent conditions holds:

$(A)$  The difference $2D$ Schr\"odinger equation
\begin{equation}\label{laxdd}
  \psi_{n+1,m+1}-u_{n,m}(\psi_{n+1,m}-\psi_{n,m+1})-\psi_{n,m}=0
\end{equation}
with
\begin{equation}\label{ud}
  u_{n,m}:=C_{nm} {\theta((n+1)U+mV+\nu_{nm}
  W+Z)\,\theta(nU+(m+1)V+\nu_{nm} W+Z)\over
  \theta((n+1)U+(m+1)V+(1-\nu_{nm})W+Z)\,\theta(nU+mV+(1-\nu_{nm})W+Z)}\,,
\end{equation}
where
\begin{equation}\label{ud11}
  2\nu_{nm}:=1+(-1)^{n+m+1}, \ \
  C_{nm}:=c_3\left(c_2^{2n+1}c_1^{2m+1}\right)^{1-2\nu_{nm}}\,
\end{equation}
and
\begin{equation}\label{pd}
  \psi_{n,m}:={\theta(A+nU+mV+\nu_{nm} W+Z)\over
  \theta(nU+mV+(1-\nu_{nm})W+Z)}\, w_1^nw_2^mw_3^{\nu_{nm}}
  \left(c_1^mc_2^n\right)^{1-2\nu_{nm}}.
\end{equation}
is satisfied for all $Z\in X$.

\medskip
$(B)$ The following identity holds:
\begin{eqnarray}\label{gr1d}
w_1w_2(c_1c_2)^{\pm 1}\wt K\left({A+U+V\mp W\over
2}\right)-w_1c_3(w_3c_1)^{\pm1}\wt K\left({A+U-V\pm W\over 2}\right)\nonumber\\
+w_2c_3(w_3c_2)^{\pm 1}\wt K\left({A+V-U\pm W \over 2}\right)-
\wt K\left({A-U-V\mp W\over 2}\right) =0\, ,\hskip16mm
\end{eqnarray}
where $\wt K:\bbC^g\to\bbC^{2^g}$ is the lifting of the Kummer map to the
universal cover.

\medskip
$(C)$ The two equations (one for the top choice of signs everywhere,
and one --- for the bottom)
\begin{eqnarray}\label{cm7d}
c_1^{\mp 2}c_3^2\ \theta(Z+U-V)\,\theta(Z-U\pm W)\,\,\theta(Z+V\pm W)&\nonumber\\
+c_2^{\mp 2}c_3^2\ \theta(Z-U+V)\,\theta(Z+U\pm W)\,\,\theta(Z-V\pm W)\nonumber\\
=c_1^{\mp 2}c_2^{\mp 2}\,\theta(Z-U-V)\,\theta(Z+U\pm W)\,\,\theta(Z+V\pm W)&\nonumber\\
+\theta(Z+U+V)\,\theta(Z-U\pm W)\,\,\theta(Z-V\pm W)& \end{eqnarray}
are valid on the theta divisor $\{Z\in X: \theta(Z)=0\}$.
\end{theo}

The purely geometric statement of part $(B)$ of this result is as
follows.
\begin{cor}[Geometric characterization of Pryms]
Given $(X,\theta)\in\A_g$, if there exist distinct points
$p_1,p_2,p_3,p_4\in
X$, none of them points of order two, such that the Kummer images of
the eight points $p_1\pm p_2\pm p_3\pm p_4$ lie on two quadrisecants
(the four points with the same parity of the number of plus signs
forming each quadruple), then $(X,\theta)$ lies in the closure of
the locus of Prym varieties of unramified double covers.
\end{cor}
\bpf Indeed, statement $(B)$ gives the two linear dependencies for
the Kummer images of the two quadruples of point. The 6 coefficients
of linear dependence appearing in these two equations depend on 6
parameters $c_i, w_i$ and are independent (since all $c_i,w_i$ can
be recovered from the 6 coefficients); thus $(B)$ says that any ppav
admitting a symmetric pair of quadrisecants is a Prym. \epf

The equivalence of $(A)$ and $(B)$ is a direct corollary of the
addition formula for the theta function. The ``only if'' part of $(A)$ is what we prove in section 2. The statement $(C)$ is actually what we
use for the proof of the ``if'' part of the theorem. The characterization of Pryms by
$(C)$ is stronger than the characterization by $(A)$. The
implication $(A)\Rightarrow (C)$ does not require the explicit
theta-functional formula for $\psi$. It is enough to require only
that equation (\ref{laxdd}) with $u$ as in (\ref{ud}) has {\it local
meromorphic} solutions which are holomorphic outside the
divisor $\theta(Un+Vm+Z)=0$ (see lemma \ref{firstlemma}).

It would be interesting to try to apply our geometric
characterization of Pryms to studying other aspects of Prym geometry
and of the geometry of the Prym locus, including the Torelli problem
for Pryms, higher-dimensional secancy conditions, representability
of homology classes in Pryms, etc. It is also tempting to ask
whether a similar characterization of Prym-Tyurin varieties of
higher order may be obtained, or whether one could use secancy
conditions to geometrically stratify the moduli space of ppavs. We
hope to pursue these questions in the future.

\section{Potential reduction of the algebro-geometric $2D$
difference Schr\"odinger operators}

To begin with let us recall a construction of algebro-geometric
difference Schr\"odinger operators proposed in \cite{kr5} (see
details in \cite{kwz}).

\noindent {\bf General notations, Baker-Akhiezer functions.}

Let $\G$ be a smooth algebraic curve of genus $\hat g$. Fix four
points $P_1^{\, \pm}, P_2^{\, \pm}\in \G$, and let $\hat
D=\g_1+\cdots+\g_{\wh g}$ be a generic effective divisor on $\G$ of
degree $\wh g$. We denote by $B$ the period matrix of the curve $\G$
(the integrals of a basis of the space of abelian differentials on
$\G$ over the $b$-cycles, once the integrals over the $a$-cycles are
normalized), by $J(\G)=\bbC^{\wh g}/\bbZ^{\wh g}+B\bbZ^{\wh g}$ --- the Jacobian
variety of $\G$, and by $\wh A:\G\hookrightarrow J(\G)$ the
Abel-Jacobi embedding of the curve into its Jacobian. We further
denote by $$\wh \theta(z):=\theta(B,z),$$ the
Riemann theta function of the variable $z\in\bbC^{\wh g}$ .

By the Riemann-Roch theorem one computes $h^0(\hat
D+n(P_1^+-P_1^-)+m(P_2^+-P_2^-))=1$, for any $n,m\in\bbZ$, and for
$\hat D$ generic. We denote by $\wh\psi_{n,m}(P),\ P\in \G$ the
unique section of this bundle. This means that $\wh\psi_{n,m}$ is the unique up to a constant factor meromorphic function such that (away from the marked points $P_i^{\pm}$) it has poles only at $\g_s$, of
multiplicity not greater than the multiplicity of $\g_s$ in $\wh D$, while at the points $P_1^+, P_2^+$ (resp. $P_1^-,P_2^-$) the function $\wh\psi_{n,m}$ has poles (resp. zeros) of orders $n$ and
$m$.

If we fix local coordinates $k^{-1}$ in the neighborhoods of
marked points (it is customary in the subject to think of marked
points as punctures, and thus it is common to use coordinates such
that $k$ at the marked point is infinite rather than zero), then
the Laurent series for $\psi_{n,m}(P)$, for $P\in\G$ near a marked
point, has the form
\begin{eqnarray}
\wh\psi_{n,m}&=&k^{\pm n}\left(\sum_{s=0}^{\infty}\xi_s^{\,
\pm}(n,m)k^{-s}\right), \ \ k=k(P), \
P\to P_1^{\, \pm}, \label{2}\\
\wh\psi_{n,m}&=&k^{\pm m}\left(\sum_{s=0}^{\infty}\chi_s^{\,
\pm}(n,m)k^{-s}\right), \ \ k=k(P), \ P\to P_2^{\, \pm}. \label{3}
\end{eqnarray}
Any meromorphic function on a Riemann surface can be expressed in
terms of the theta functions, but it is easier to write an
expression for $\wh\psi_{n,m}$ using both the theta functions and
the differentials of the third kind. Indeed, for $i=1,2$ let
$d\wh\Omega^i\in H^0(K_\G+P_i^++P_i^-)$ be the differential of the
third kind, normalized to have residues $\mp 1$ at $P_i^{\, \pm}$ and with zero integrals over all the $a$-cycles,
and let $\wh\Omega^i$ be the corresponding abelian integral, i.e.
the function on the Riemann surface obtained by integrating $d\wh
\Omega^i$ from some fixed starting point to the variable point. Then
we have the following expression
\begin{equation}\label{4}
\wh\psi_{n,m}(P)=r_{nm}{\wh\theta(\wh A(P)+n\wh U+m\wh V+\wh Z)
\over \wh \theta(\wh A(P)+\wh Z)\,}\
e^{n\wh\Omega_1(P)+m\wh\Omega_2(P)},
\end{equation}
where $r_{nm}$ is some constant, $\wh U=\wh A(P_1^-)-\wh A(P_1^+), \
\ \wh V=\wh A(P_2^-)-\wh A(P_2^+),$ and
\begin{equation}\label{hatz}
\wh Z=-\sum_s \wh A(\g_s)+\wh\kappa,
\end{equation}
where $\wh\kappa$ is the vector of Riemann constants. Indeed, to
prove that such an expression for $\wh\psi_{n,m}$ is valid, one
only needs to verify that both sides have the same zeros and poles,
which is clear by construction.

{\bf Notation.} From now on it will be useful to think of $n$ and
$m$ as discrete variables, which are shifted by the shift operators
that we denote $T_1:n\mapsto n+1$ and $T_2:m\mapsto m+1$
respectively. To emphasize the difference between the operator and
its action, for a function $f=f(n,m)$ we will write ${\bf t}_\mu
f:=T_\mu\circ f$, so that for example $T_1(f\cdot g)={\bf t}_1f\cdot
{\bf t}_1g$. We will also denote $H:=T_1T_2-u(T_1-T_2)-1$ the
difference operator that is very important for what follows.

\begin{theo}[\cite{kr5}]
The Baker-Akhiezer function $\wh \psi_{n,m}$ given by formula (\ref{4})
satisfies the following difference equation
\begin{equation}\label{eqn}
\wh\psi_{n+1,m+1}-a_{n,m}\wh\psi_{n+1,m}-b_{n,m}\wh\psi_{n,m+1}
+c_{n,m}\wh\psi_{n,m}=0,
\end{equation}
where we let
\begin{equation}\label{coef}
a_{n,m}:={\xi_0^+(n+1,m+1)\over \xi_0^+(n+1,m)}\,, \ \
b_{n,m}:={\chi_0^+(n+1,m+1)\over \chi_0^+(n,m+1)},
\end{equation}
\begin{equation}\label{coef1}
c_{n,m}:=b_{n,m}{\xi^-(n,m+1)\over\xi_0^-(n,m)}=
{\xi^-(n,m+1)\,\chi_0^+(n+1,m+1)\over \xi_0^-(n,m)\,\chi_0^+(n,m+1)}.
\end{equation}
\end{theo}
Explicit $\theta$-functional formulae for the coefficients follow
from equation (\ref{4}) which implies
\begin{equation}\label{coef2a}
  \xi_0^{\, \pm}=r_{nm}{\wh\theta(\wh A(P_1^{\, \pm})+n\wh U+m\wh V+\wh
  Z)\, \over \wh\theta(\wh A(P_1^{\, \pm})+\wh Z)}\ e^{n\a_{1}^{\,
  \pm}+m\a_2^{\, \pm}}
\end{equation}
\begin{equation}\label{coef2}
  \chi_0^{\, \pm}=r_{nm}{\wh\theta(\wh A(P_2^{\, \pm})+n\wh U+m\wh V+\wh
  Z)\over\wh\theta(\wh A(P_2^{\, \pm})+\wh Z)} e^{n\b_{1}^{\,
  \pm}+m\b_2^{\, \pm}}
\end{equation}
The constants $\a_i^{\, \pm},\b_i^{\, \pm}$ are defined by the formulae:
\begin{eqnarray}
  \a_2^{\pm}=\Omega_2(P_1^{\,\pm}); \ \ \Omega_1&=\pm\ln k+\a_1^{\, \pm}
  +O(k^{-1}),\ \ P\to P_1^{\, \pm}, \label{coef3}\\
  \b_1^{\pm}=\Omega_1(P_2^{\,\pm}); \ \ \Omega_2&=\pm\ln k+\b_2^{\, \pm}
  +O(k^{-1}),\ \ P\to P_2^{\, \pm}. \label{coef4}
\end{eqnarray}

\bigskip
\noindent {\bf Setup for the Prym construction}

We now assume that the curve $\G$ is an algebraic curve endowed with
an involution $\s$ without fixed points; then $\G$ is a unramified
double cover $\G\longmapsto \G_0$, where $\G_0=\G/\s$. If $\G$ is of
genus $\wh g=2g+1$, then by Riemann-Hurwitz the genus of $\G_0$ is
$g+1$. On $\G$ one can choose a basis of cycles $a_i, b_i$ with the
canonical matrix of intersections $a_i\cdot a_j=b_i\cdot b_j=0,\
a_i\cdot b_j=\delta_{ij},\ \ 0\leq i,j\leq 2g,$ such that under the
involution $\s$ we have $\s(a_0)=a_0,\ \s(b_0)=b_0,\
\s(a_j)=a_{g+j},\ \s(b_j)=b_{g+j}, 1\leq j\leq g$. If $d\omega_i$
are normalized holomorphic differentials on $\G$ dual to this choice
of $a$-cycles, then the differentials
$du_j=d\omega_j-d\omega_{g+j}$, for $j=1\ldots g$ are odd, i.e.
satisfy $\s^*(du_k)=-du_k$, and we call them the normalized
holomorphic Prym differentials. The matrix of their $b$-periods
\begin{equation}\label{pi}
  \Pi_{kj}=\oint_{b_k}du_j,\ \ 1\leq k,j\leq g\,,
\end{equation}
is symmetric, has positive definite imaginary part, and defines the
Prym variety $$\P(\G):=\bbC^g/\bbZ^g+\Pi\bbZ^g$$ and the
corresponding Prym theta function $$\theta(z):=\theta(\Pi,z),$$ for
$z\in\bbC^g$. We assume that the marked points $P_1^{\, \pm},
P_2^{\, \pm}$ on $\G$ are permuted by the involution, i.e.
$P_i^+=\s(P_i^-)$. For further use let us fix in addition a third
pair of points $P_3^{\pm},$ such that also $P_3^-=\s(\P_3^+)$.

The Abel-Jacobi map $\G\hookrightarrow J(\G)$ induces the Abel-Prym
map $A:\G\longmapsto\P(\G)$. There is a choice of the base point
involved in defining the Abel-Jacobi map, and thus in the Abel-Prym
map; let us choose this base point (such a choice is unique up to a
point of order two in $\P(\G)$) in such a way that
\begin{equation}\label{prab}
  A(P)=-A(\s(P)).
\end{equation}

\smallskip
\noindent {\bf Admissible divisors.} An effective divisor on $\G$
of degree $\hat g-1=2g$, $D=\g_1+\ldots\g_{2g}$, is called {\it
admissible} if it satisfies
\begin{equation}\label{const}
  [D]+[\s(D)]=K_\G\in J(\G)
\end{equation}
(where $K_\G$ is the canonical class of $\G$), and if moreover
$H^0(D+\s(D))$ is generated by an even holomorphic differential
$d\Omega$, i.e. that
\begin{equation}\label{const0}
  d\Omega(\g_s)=d\Omega(\s(\g_s))=0,\ \  d\Omega=\s(d\Omega).
\end{equation}

Algebraically, what we are saying is the following. The divisors $D$
satisfying (\ref{const}) are the preimage of the point $K_\G$ under
the map $1+\s$, and thus are a translate of the subgroup
$Ker(1+\s)\subset J(\G)$ by some vector. As shown by Mumford
\cite{mum2}, this kernel has two components --- one of them being
the Prym, and the other being the translate of the Prym variety by
the point of order two corresponding to the cover $\Gamma\to
\Gamma_0$ as an element in $\pi_1(\Gamma_0)$. The existence of an
even differential as above picks out one of the two components, and
the other one is obtained by adding $A-\sigma(A)$ to the divisor of
such a differential, for some $A$. In Mumford's notations the
component we pick is in fact $P^-$ (when we choose the base point according to (\ref{prab}) to identify ${\rm Pic}^0$ and ${\rm Pic}^{\wh g -1}$), but throughout this paper we
will have to deal with both components, using some point (which will
be called $P_3^+$) and the corresponding shift by $P_3^+-P_3^-$ to
pass from one component to the other. We will prove the following
statement.

\begin{prop}\label{theta.adm}
For a generic vector $Z$ the zero-divisor $D$ of the function
$\theta(A(P)+Z)$ on $\G$ is of degree $2g$ and satisfies the
constraints (\ref{const}) and (\ref{const0}), i.e. is admissible.
\end{prop}
{\bf Remark.} We have been unable to find a proof of this statement
in the literature. However, both Elham Izadi and Roy Smith have
independently supplied us with simple proofs of this result, based
on Mumford's description and results on Prym varieties. The reason
we have chosen to still give the longer analytic proof below is
because we need some of the intermediate results later on, and
also to give an independent analytic proof of some of Mumford's
results.

Note that the function $\theta(A(P)+Z)$ is multi-valued on $\G$, but
its zero-divisor is well-defined. The arguments identical to that in
the standard proof of the inversion formula (\ref{hatz}) show that
the zero divisor $D(Z):=\theta(A(P)+Z)$ is of degree $\hat g-1=2g$.

\begin{lem} For a generic $D=D(Z)$ and for each set of integers $(n,m,r)$
such that
\begin{equation}\label{ev}
n+m+r=0  \mod 2
\end{equation}
the space
$$
  H^0(D+n(P_1^+-P_1^-)+m(P_2^+-P_2^-)+r(P_3^+-P_3^-))
$$
is one-dimensional. A basis element of this space is given by
\begin{equation}\label{psiprym}
\psi_{n,m,r}(P):=h_{n,m,r}{\theta(A(P)+nU+mV+rW+ Z)\over\theta( A(P)+Z)}\
e^{n\Omega_1(P)+m\Omega_2(P)+r\Omega_3(P)},
\end{equation}
where $\Omega_j$ is the abelian integral corresponding to a differential
$d\Omega_i$ of the third
kind, odd under the involution $\s$, and with residues $\mp 1$ at $P_j^{\pm}$ (i.e.
$d\Omega_i=-\s(d\Omega_j)$), satisfying the normalization condition
\begin{equation}\label{a-cycle}
\oint_{a_k}d\Omega_j=\pi i \,l_k,\ \ l_k\in \bbZ, \ \ k=0,\ldots, 2g,
\end{equation}
and $U,V,W$ are the vectors of $b$-periods of these differentials,
i.e.
\begin{equation}\label{per}
2\pi i U_k=\oint_{b_k}d\Omega_1,\
2\pi i V_k=\oint_{b_k}d\Omega_2,\
2\pi i W_k=\oint_{b_k}d\Omega_3.
\end{equation}
\end{lem}
\bpf It is easy to check that the right hand side of (\ref{psiprym})
is a single valued function on $\G$ having all the desired
properties, and thus it gives a section of the desired bundle. Note
that the constraint (\ref{ev}) is required due to (\ref{a-cycle}),
and the uniqueness of $\psi$ up to a constant factor, i.e. the
one-dimensionality of the $H^0$ above, is a direct corollary of the
Riemann-Roch theorem. \epf

For further use let us note that bilinear Riemann identities imply
\begin{equation}\label{periods}
  2U=A(P_1^-)-A(P_1^+), \  \ 2V=A(P_2^-)-A(P_2^+), \ \ 2W=A(P_3^-)-A(P_3^+).
\end{equation}

Let us compare the definition of $\wh \psi_{n,m}$ defined for any
curve $\G$, with that of $\psi_{n,m,r}$, which is only defined for a
curve with an involution satisfying a number of conditions. To make
such a comparison, consider the divisor $\wh D=D+P_3^+$ of degree
$\hat g=2g+1$, and let $\wh\psi_{n,m}$ be the corresponding
Baker-Akhiezer function.
\begin{cor}
For the Baker-Akhiezer function $\wh \psi_{nm}$ corresponding to the
divisor $\wh D=D+P_3^+$ we have
\begin{equation}\label{psi-psi}
  \wh \psi_{nm}=\psi_{n,m,\nu}
\end{equation}
where $\nu=\nu_{nm}$ is defined in (\ref{ud11}), i.e. is 0 or 1 so that
$n+m+\nu$ is even.
\end{cor}
\begin{cor}
If $n+m$ is even, then by formulae (\ref{4},\ref{psiprym})
\begin{eqnarray}\label{comp}
 {\wh\theta(\wh A(P)+n\wh U+m\wh V+\wh Z)
  \,\wh \theta(\wh A(P_0)+\wh Z)\over
  \wh \theta(\wh A(P)+\wh Z)\,
  \wh\theta(\wh A(P_0)+n\wh U+m\wh V+\wh Z)}= \ \ \ \ \ \ \   \nonumber\\
  {\theta(A(P)+nU+mV+ Z)\, \theta( A(P_0)+Z)\over
  \theta( A(P)+Z)\, \theta(A(P_0)+nU+mV+ Z)}e^{nr_1+mr_2},
\end{eqnarray}
where $r_i=\int_{P_0}^P(d\wh \Omega_i-d\Omega_i)$, and we recall
that $\wh Z=\wh A(\wh D)+\wh \kappa$, and $Z$ is its image.
\end{cor}
\begin{rem}
This equality, valid for any pair of points $P,P_0$ is a non-trivial
identity between theta functions. The authors' attempts to derive it
directly from the Schottky-Jung relations have failed so far.
\end{rem}

\smallskip
\noindent{\bf Notation.}
For brevity throughout the rest of the paper we use the notation:
$\psi_{n,m}:=\psi_{n,m,\,\nu_{nm}}$.

\begin{lem}
The Baker-Akhiezer function $\psi_{n,\,m}$ given by
\begin{equation}\label{psipr}
\psi_{n,\,m}={\theta(A(P)+Un+Vm+\nu_{nm} W+Z)\over
  \theta(Un+Vm+(1-\nu_{nm}) W + Z)\, \theta( A(P)+Z)}
  \cdot{e^{n\Omega_1(P)+m\Omega_2(P)+\nu_{nm}\Omega_3(P)}\over
 e^{(2\nu_{nm}-1)(n\Omega_1(P_3^+)+m\Omega_2(P_3^+))}},
\end{equation}
satisfies the equation (\ref{laxdd}), i.e.
$$
  \psi_{n+1,\,m+1}-u_{n,m}(\psi_{n+1,m}-\psi_{n,\,m+1})-\psi_{n,m}=0,
$$
with $u_{n,m}$ as in (\ref{ud},\ref{ud11}), where
\begin{equation}\label{d2}
  c_1=e^{\Omega_2(P_3^+)},\ \  c_2=e^{\Omega_1(P_3^+)},\ \ c_3=e^{\Omega_1(P_2^+)}
\end{equation}
\end{lem}
\bpf Note that the first and the last factors in the denominator of
(\ref{psipr}) correspond to a special choice of the normalization
constants $h_{n,m,\, \nu}$ in (\ref{psiprym}):
\begin{eqnarray}\label{nordec}
  \psi_{nm}(P_3^-)=&(\theta(Z+W))^{-1}, \ \ \nu_{nm}=0,\nonumber\\
  \psi_{nm}e^{-\Omega_3}|_{\,P=P_3^+}=&(\theta(Z-W))^{-1} ,\ \ \nu_{nm}=1.
\end{eqnarray}
This normalization implies that for even $n+m$ the difference
$(\psi_{n+1,m+1}-\psi_{n,m})$ equals zero at $P_3^-$. At the same
time as a corollary of the normalization we get that
$(\psi_{n+1,m}-\psi_{n,m+1})$ has no pole at $P_3^+$. Hence, these
two differences have the same analytical properties on $\G$ and thus
are proportional to each other (the relevant $H^0$ is
one-dimensional by Riemann-Roch). The coefficient of proportionality
$u_{nm}$ can be found by comparing the singularities of the two
functions at $P_1^+$. \epf

The second factor in the denominator of the formula (\ref{psipr})
does not affect equation (\ref{laxdd}). Hence, the lemma proves the
``only if'' part of the statement $(A)$ of the main theorem for the case
of smooth curves. It remains valid under degenerations to singular
curves which are smooth outside of fixed points $Q_k$ which are
simple double points, i.e. to the curves of type $\{\G,\s,Q_k\}$.

\begin{rem}
Equation (\ref{laxdd}) as a special reduction of (\ref{eqn}) was
introduced in \cite{grin}. It was shown that equation (\ref{eqn})
implies a five-term equation
\begin{equation}\label{51}
  \psi_{n+1,m+1}-\tilde a_{nm}\psi_{n+1,m-1}-\tilde b_{n,m}\psi_{n-1,m+1}+
  \tilde c_{nm} \psi_{n-1,m-1}
  =\tilde d_{n,m}\psi_{n,m}
\end{equation}
if and only if it is of the form (\ref{laxdd}). A reduction of the
algebro-geometric construction proposed in \cite{kr5} in the case of
algebraic curves with involution having two fixed points was found.
It was shown that the corresponding Baker-Akhiezer functions do
satisfy an equation of the form (\ref{laxdd}). Explicit formulae for
the coefficients of the equations in terms of Riemann
theta-functions were obtained. The fact that the Baker-Akhiezer
functions and the coefficients of the equations can be expressed in
terms of Prym theta-functions is new.
\end{rem}

We are now ready to complete the proof of proposition
(\ref{theta.adm}). Let $\psi_{n,m}$ be the Baker-Akhiezer function
given  by (\ref{psipr}). According to Lemma 2.3 it satisfies
equation (\ref{laxdd}). The differential $d\psi_{n,m}$ is also a
solution of the same equation, and thus we get, using the shift
operator notation,
\begin{equation}\label{nnov31}
  (T_1-1)(\psi_{n,m}^{\s}d\psi_{n,m+1}-\psi_{n,m+1}^{\s}d\psi_{n,m})=
  (T_2-1)(\psi_{n,m}^{\s}d\psi_{n+1,m}-\psi_{n+1,m}^{\s}d\psi_{n,m})
\end{equation}

For a generic set of algebro-geometrical spectral data the products
$\psi_{n,m}^{\s}\psi_{n,m+1}$ and $\psi_{n,m}^{\s}\psi_{n+1,m}$ are
quasi-periodic functions of the variables $n$ and $m$. The data for
which they are periodic are characterized as follows.

Let $dp_j,\, i=1,2$ be the third kind abelian differentials with
residues $\mp 1$ at the punctures $P_j^{\pm}$, respectively, and
normalized by the condition that all of their periods are {\it purely imaginary},
\begin{equation}\label{perimag}
  \Re\oint_c dp_j=0, \ \ \ \forall c\in H^1(\G,Z).
\end{equation}
Non-degeneracy of the imaginary part of the period matrix of
holomorphic differential implies that such $dp_j$ exists and is
unique. If the periods of $dp_j$ are of the form
\begin{equation}\label{perimag1}
  \oint_c dp_j={\pi i n_c^j\over N_j},\  n_c^j\in \bbZ,
\end{equation}
then the function $\mu_j(Q)=e^{N_j\int^Q dp_j}$ is single-valued on
$\G$, has pole of order $N_j$ at $P_j^+$ and zero of order $N_j$ at
$P_j^-$. From the uniqueness of the Baker-Akhiezer function it then
follows that
\begin{eqnarray}\label{muper}
  \psi_{n+2N_1,m}&=&{\mu_1\over \mu_1(P_3^-)}\ \psi_{n,m},\ \ \
  \psi_{n,m+2N_2}={\mu_2\over \mu_2(P_3^-)}\ \psi_{n,m}\ \ ,\ \ \nu=0 \nonumber\\
  \psi_{n+2N_1,m}&=&{\mu_1\over \mu_1(P_3^+)}\ \psi_{n,m},\ \ \
  \psi_{n,m+2N_2}={\mu_2\over \mu_2(P_3^+)}\ \psi_{n,m}\ \ ,\ \ \nu=1
\end{eqnarray}
These imply
\begin{equation}\label{nn1}
  \psi_{n+2N_1,m}^{\s}d\psi_{n+1+2N_1,m}=
  \psi_{n,m}^{\s}d\psi_{n+1,m}+(\psi_{n,m}^{\s}\psi_{n+1,m})dp_1
\end{equation}
and similar monodromy properties for the other terms in
(\ref{nnov31}). In this case the averaging of equation
(\ref{nnov31}) in the variables $n,m$ gives the equation
\begin{equation}\label{nnov32}
  \langle\psi^{\s}{\bf t}_2\psi-{\bf t}_2\psi^{\s})\psi\rangle_2dp_1=
  \langle\psi^{\s}({\bf t}_1\psi)-({\bf t}_1\psi^{\s})\psi\rangle_1\, dp_2.
\end{equation}
Here $\langle\cdot\rangle_1$ stands for the mean value in $n$
and $\langle\cdot\rangle_2$ stands for the mean value in $m$.
For a generic curve  differentials $dp_j$ have no common zeros.
Hence, for such curves the differential
\begin{equation}\label{nnov33}
  d\Omega={dp_1\over
  \langle\psi^{\s}{\bf t}_1\psi-{\bf t}_1\psi^{\s}\psi\rangle_1}=
  {dp_2\over\langle\psi^{\s}{\bf t}_2\psi-{\bf t}_2\psi^{\s}\psi\rangle_2}
\end{equation}
is holomorphic on $\G$. It has zeros at the poles of $\psi$ and
$\psi^{\s}$. The curves for which (\ref{perimag1}) holds for some
$N_j$ are dense in the moduli space of all smooth genus $g$ curves.
That proves that equation (\ref{nnov33}) holds for any curve.
Proposition (\ref{theta.adm}) is proven.
\begin{rem}
We have thus proven that for any Prym variety part (A) of the main theorem is satisfied. Note, however, that the statement of the main theorem is for all abelian varieties in the closure of the locus $\P_g$ in $\A_g$. To show that condition (A) holds for abelian varieties in the closure, it is enough to note that (A) is an algebraic condition, and thus is valid on the closure of the locus.

When we prove the characterization --- the ``only if'' part of the main theorem --- in section 5, there will be no problems with the closure as we will be able to show explicitly that condition (C) (implied by (A)) exhibits the abelian variety as the Prym for a possibly nodal curve.
\end{rem}

\section{A discrete analog of Novikov-Veselov hierarchy}
In this section we introduce {\it multi-parametric} deformations of
the Baker-Akhiezer functions and prove that they satisfy a system of
difference-differential equations. The compatibility conditions of
these equations can be regarded as a discrete analog of the
Novikov-Veselov hierarchy (\cite{nv}).

Let $t=\{t_i^1,t_i^2,\ i=1,2\ldots\}$ be two sequences of complex
numbers (we assume that only finitely many of them are non-zero).
We will construct a function $\psi$ on the curve $\G$ with prescribed
exponential essential singularities at the points $P_i^\pm$
controlled by these $t$.

\begin{lem}
Let $D=D(Z)=\gamma_1+\ldots+\gamma_{2g}$ be an admissible divisor.
Then there exists a unique up to a constant factor meromorphic
function $\psi_{n,m}(t,P)$  of $P\in\Gamma$, which we call a
multi-parametric deformation of the Baker-Akhiezer function, such
that
\begin{itemize}
\item[(i)] outside of the marked points it has poles only at the points
 $\g_s$ of multiplicity not greater than the multiplicity of
 $\g_s$ in $D$
\item[(ii)] $\psi_{n,m}(t,P)$ has an at most simple pole at $P_3^+$
\item[(iii)] in the local coordinate $k^{-1}$ mapping a small neighborhood of
 $P_1^{\pm}$ to a small disk in $\bbC$, (with the marked point
 mapping to zero),  it has the power series expansion
\begin{equation}\label{200}
  \psi_{n,m}(t,P)=k^{\mp n}e^{\pm\sum\limits_i
  t_i^1k^{-i}}\left(\sum_{s=0}^{\infty}\xi_s^{\,
  \pm}(n,m,t)k^s\right),
\end{equation}
for some $\xi_s^\pm$ (notice that this means there is an
essential singularity, and the expansion starts from $k^{-n}$ at
$P_1^+$ and $k^n$ at $P_1^-$, and goes towards $k^{-\infty}$)
\item[(iv)] in the local coordinate $k^{-1}$ near $P_2^{\pm}$ it has the
 power series expansion
\begin{equation}\label{300}
  \psi_{n,m}(t,P)=k^{\mp m}e^{\pm\sum\limits_i
  t_i^2k^{-i}}\left(\sum_{s=0}^{\infty}\chi_s^{\,
  \pm}(n,m,t)k^s\right).
\end{equation}
\end{itemize}
\end{lem}
\bpf
This function $\psi_{n,m}$ is given by
\begin{eqnarray}\label{psiprym10}
  \psi_{n,m}(t,P)=h_{n,m}(t){\theta(A(P)+nU+mV+\nu_{nm}W+Z+\sum_i(t_i^1 U_i^1+t_i^2 U_i^2))
  \over\theta( A(P)+Z)}\times\nonumber\\
  \times\exp\left(n\Omega_1(P)+m\Omega_2(P)+\nu_{nm}\Omega_3(P)+
  \sum_i (t_i^1 \Omega_i^1(P)+t_i^2 \Omega_i^2(P))\right),
\end{eqnarray}
where $\Omega_1,\Omega_2,\Omega_3$ and the vectors $U,V,W$ are as in
Lemma 2.2; $\Omega_j^\mu$ for $\mu=1,2$ is the abelian integral of
the differential $d\Omega_j^{\mu}$ which has poles of the form
\begin{equation}\label{dec1}
  d\Omega_j^{1(2)}=\pm d(k^j+O(1))
\end{equation}
at the punctures $P_{1(2)}^{\pm}$, holomorphic everywhere else and
is uniquely determined by the normalization conditions
\begin{equation}\label{a-cycle1}
  \oint_{a_k}d\Omega_j^{1(2)}=0, \ \ k=0,\ldots, 2g\,;
\end{equation}
coordinates of the vectors $U_j^{1(2)}$ are defined by $b$-periods
of these differentials, i.e.
\begin{equation}\label{per12}
  2\pi i U_{k,j}^{1(2)}=\oint_{b_k}d\Omega_j^{1(2)},\ k=1,\ldots,g
\end{equation}
Note that, as before, if $\nu_{nm}=0$ then $\psi_{n,m}$ is in fact
holomorphic at $P_3^+$, and if $\nu_{nm}=1$, then $\psi_{n,m}$ does
have a pole at $P_3^+$, but also has a zero at $P_3^-$. As before,
we normalize $\psi_{n,m}$ by the conditions (\ref{nordec}). \epf

\noindent{\bf Notations.} In what follows we will deal with formal
pseudodifference operators, shifting $n$ and $m$, with
coefficients being functions of the variables $n$ and $m$, and of
the $t$'s. From now on when we write functions $f,g,\ldots$ as
coefficients of pseudodifference operators, they are meant to be
functions of $n,m$ and $t$.

Denote by ${\mathcal R}$ the ring of functions of variables $n,m,$
and $t$. We denote by $\O_1^{\pm}$ the rings of pseudodifference
operators in two variables that are Laurent polynomials in
$T_1^{\pm}$, i.e.
$$
  \O_1^{\pm}:={\mathcal R}((T_1^{\mp}))[T_2,T_2^{-1}]=\lbrace D=\sum_{j=M_1}^{M_2}\sum_{i=N}^\infty r_{ij}T_1^i T_2^j\rbrace,
$$
where $r_{ij}\in{\mathcal R}$. The intersection
$$\O:=\O_1^+\cap\O_1^-={\mathcal R}[T_1,T_1^{-1},T_2,T_2^{-1}]$$ is
the ring of difference operators. We further denote $\O_{1,0}^\pm$
the ring of pseudodifference operators in one variable that are
Laurent polynomials in $T_1^\mp$, thought of as subrings of
$\O_1^\pm$, respectively, i.e.
$$
  \O_{1,0}^\pm:={\mathcal R}((T_1^\mp))=\lbrace D=\sum_{i=N}^\infty
r_{i}T_1^i\rbrace.
$$
Finally we denote by $\O_H$ the left principal ideal generated by
the operator $H=T_1T_2-u(T_1-T_2)-1$, i.e. $\O_H:=\O H$, and
similarly set $\O_H^\pm:=\O_{1}^\pm H$.

Moreover, in doing computations in these rings it is often
convenient to compute only a couple highest terms. To this end, we
will use for $k>0$ notations $O(T_1^{-k})=T_1^{-k}{\mathcal
R}[[T_1^{-1}]]$ for the operators in $\O_1^+$ only having the terms
with $T_1^n$ for $n\le -k$, and by $O(T_1^k)=T_1^k{\mathcal
R}[[T_1]]$ --- operators in $\O_1^-$ only having terms with $T_1^n$
for $n\ge k$.

\medskip

We now want to show that the multi-parametric deformations of the
Baker-Akhiezer functions satisfy a hierarchy of
difference-differential equations --- the result is as follows.
\begin{prop}\label{lmL}
The Baker-Akhiezer function $\psi=\psi_{n,m}(t,P)$ satisfies
(\ref{laxdd}) with  $u_{nm}$ as in (\ref{ud}), with $Z$ replaced by
$Z+\sum_j(t_j^1 U_j^1+t_j^2 U_j^2)$ (this can be written as
$H\psi=0$). There exist unique difference operators of the form
\begin{equation}\label{Li}
  L_{j}^{(\k)}=\left(f_{0j}+\sum_{i=1}^{j-1} f_{ij}^{(\k)}
  T_{\k}^{\,i}+T_{\k}^{-i}f_{ij}^{(\k)}\right)
  \left(T_{\k}-T_{\k}^{-1}\right),\ \ \k=1,2,\ \ j=1,2,\ldots
\end{equation}
such that the equations
\begin{equation}\label{dec25}
  {\p\over\p t_j^{\k}}\ \psi=L_j^{(\k)}\psi
\end{equation}
hold.
\end{prop}
\bpf
The proof of the first statement is identical to that in lemma 2.3. The proof
of the statement that there are operators of the form
\begin{equation}\label{Li1}
  L_{j}^{(\k)}=
  \sum_{i=-j}^{j} g_{ij}^{(\k)}T_{\k}^{\,i}
\end{equation}
such that (\ref{dec25}) hold is standard. Indeed, for each formal series (\ref{200})
there exists a unique operator $L_{j}^{(1)}$ such that
\begin{equation}
  \left({\p\over\p t_j^{\k}}-L_j^{(1)}\right)\psi=k^{\pm n}e^{\pm\sum_i t_i^1k^i}
  \left(\sum_{s=0}^{\infty}\tilde\xi_s^{\,\pm}k^{-s}\right),\ \ \tilde\xi_0^+=0.\label{201}
\end{equation}
The coefficients $g_j^{(1)}$ of the operator are difference
polynomials in terms of  the coefficients $\xi_s$ of the series
(\ref{200}). Now note that the left-hand-side of (\ref{201})
satisfies all the properties that $\psi$ satisfies, and thus must be
proportional to it. However, since $\tilde\xi_0^+=0$, the constant
of proportionality must be equal to zero, and thus the
left-hand-side vanishes as desired. The same arguments proves the
existence of $L_j^{(2)}$.

It remains to show that the operators $L_j^{(\k)}$ have the form
(\ref{Li}). This is a matter of showing that the coefficients of
$L_j^{(\k)}$ satisfy certain identities, i.e. that the generally
constructed $g_{ij}^{(\k)}$ for $-j\le i\le j$ can in fact be
expressed in terms of $f_{ij}^{(\k)}$ for $1\le i\le j-1$. One
easily checks that if $L_j^{(\k)}$ is given by (\ref{Li}), then
its operator formal adjoint satisfies
\begin{equation}\label{adj}
  \left(L_j^{(\k)}\right)^*=-(T_{\k}-T_{\k}^{-1})L_j^{(\k)}(T_{\k}-T_{\k}^{-1})^{-1}.
\end{equation}
This equation is in fact equivalent to (\ref{Li}), as it determines
the coefficients of all the negative powers of $T_\k$ uniquely,
given the coefficients of the positive powers. It thus remains to
prove this identity.

\smallskip
We denote by $\psi^\s$ the composition $\psi^\s(P)=\psi(\s(P))$
--- notice that by (\ref{200}) we know the expansion of both
$\psi$ and $\psi^\s$ near $P_1^\pm$. Now consider the differential
$\psi^{\s}\left(T_{\k}^n(T_{\k}-T_{\k}^{-1})\,\psi\right)d\Omega$,
where, as before, $d\Omega$ is a holomorphic differential having
zeros at poles of $\psi$ and $\psi^{\s}$. Then this expression is
a meromorphic differential on $\G$ which a priori has poles only
at $P_{\k}^{\pm}$ and $P_3^{\pm}$. Due to normalization
(\ref{nordec}) it is holomorphic at the punctures $P_3^{\pm}$ ---
the pole of $\psi$ cancels with the zero of $\psi^\s$ and vice
versa. Therefore, for $n>0$ this differential has a pole only at
$P_{\k}^+$, and hence its residue at this point must vanish:
\begin{equation}\label{dec311}
  \res_{P_{\k}^+}\left(\psi^{\s}\left(T_{\k}^n(T_{\k}-T_{\k}^{-1})\,
  \psi\right)d\Omega\right)=0, \quad \forall n>0.
\end{equation}
The normalization (\ref{nordec}) implies also
\begin{equation}\label{dec312}
  \res_{P_{\k}^+}\left(\psi^{\s}\left(T_{\k}\,\psi\right)d\Omega\right)=1.
\end{equation}
Equations (\ref{dec311},\ref{dec312}) recurrently define
coefficients of the power series expansion of $\psi^{\s}d\Omega$ at
$P_{\mu}^+$ in terms of the coefficients of the power series for
$\psi$. The corresponding expressions can be explicitly written in
terms of the so-called wave operator.

We first observe that in the ring $\O_{1,0}^{+}$ there exists a
unique pseudo-difference operator
\begin{equation}\label{decphi}
\Phi=\sum_{s=0}^{\infty}\varphi_s T_1^{-s}
\end{equation}
such that the expansion (\ref{200}) of $\psi$ at $P_1^+$ is equal to
\begin{equation}\label{dec313}
  \psi=\Phi \ k^n e^{\sum_i t_i^1k^i}.
\end{equation}
Indeed, this identity gives a unique way to determine the coefficients $\varphi_s$ recursively.
\begin{lem}
The following identity holds:
\begin{equation}\label{dec314}
  \psi^{\s}d\Omega=\left(k^{-n} e^{-\sum_i t_i^1k^i} (T_1-T_1^{-1})\ \Phi^{-1}
  (T_1-T_1^{-1})^{-1} \right){dk\over k^2-1}
\end{equation}
\end{lem}
Here and below the right action of pseudo-difference operators is
defined as the formal adjoint action, i.e we set $fT=T^{-1}f$.

\bpf Recall that by definition the residue of a pseudo-differential
operator $D=\sum_s d_sT^s$ is $\res_T D:=d_0$. It is easy to check
--- by verifying that this holds for the basis, i.e. checking this
for $D_1=T_1^a$ and $D_2=T_1^b$ --- that for any two
pseudo-differential operators $D_1,D_2$ we have
\begin{equation}\label{dec315}
\res_k \left(k^{-n} e^{-\sum_i t_i^1k^i} D_1\right) \left(D_2 k^n
e^{\sum_i t_i^1k^i}\right)\,{d\ln k}=\res_T \left(D_2D_1\right)
\end{equation}
The last equation implies that
$$\res_k\left(k^{-n} e^{-\sum_i t_i^1k^i} (T_1-T_1^{-1})\ \Phi^{-1}
  (T_1-T_1^{-1})^{-1} \right)\left(T_1^n(T_1-T_1^{-1})\,\psi\right){dk\over k^2-1}=
$$
$$
\res_k\left(k^{-n} e^{-\sum_i t_i^1k^i} \Phi^{-1}
  (T_1-T_1^{-1})^{-1} \right)\left(T_1^n(T_1-T_1^{-1})\Phi \,k^n
e^{\sum_i t_i^1k^i}\right)d\ln k=\res_T\  T_1^n=\delta_{n,\,0},
$$
i.e. the formal series defined by the right hand side of
(\ref{dec314}) satisfies the equations
(\ref{dec311},\ref{dec312}), which are the defining equations (by
solving term by term, see above) for $\psi^{\s}d\Omega$. \epf

Now we are ready to complete the proof that the adjoints of
$L_j^{\mu}$ satisfy (\ref{adj}), thus proving proposition \ref{lmL}.
Consider the pseudo-difference operator
$$
 \L:=\Phi T_1 \Phi^{-1},
$$
for which $\psi$ is an eigenvector: indeed
\begin{equation}\label{dec316}
  \L\psi=\Phi T_1 k^n e^{\sum_i t_i^1k^i}=\Phi k^{n+1} e^{\sum_i
t_i^1k^i}=k\psi.
\end{equation}
Considering the expansion of (\ref{dec25}) in a neighborhood of $P_1^+$, we see that
the positive parts of the pseudo-difference operators $L_j^{(1)}$ and $\L^j$
coincide:
\begin{equation}\label{d+}
  (L_j^{(1)})_+=\L^j_+
\end{equation}
(where by the positive part of a pseudo-difference operator
$D=\sum_s d_sT^s$ we mean $D_+:=\sum_{s>0} d_sT^s$).

The differential $d\Omega$ is independent of $n$. Therefore, from
(\ref{dec314}) it follows that the operator
\begin{equation}
  \wt\L:=(T_1-T_1^{-1})^{-1}\L^*(T_1-T_1^{-1}).
\end{equation}
has $\psi^\s$ as an eigenfunction:
\begin{equation}
  \wt\L\psi^{\s}=k\psi^{\s}.
\end{equation}
Equation (\ref{dec25}) considered in the neighborhood of $P_1^-$ implies
that the negative parts of $L_j^{(1)}$ and $\wt \L^j$ coincide,
\begin{equation}\label{dec318}
  (L_j^{(1)})_-= -\wt \L^j_-
\end{equation}
The last two equations prove (\ref{adj}) and then (\ref{Li}) for
$\k=1$. The case $\k=2$ is analogous, and the proposition is thus proven.
\epf

\begin{cor}
The operators $H$ and $L_{j}^{(\k)}$ satisfy the equations
\begin{equation}\label{nvd}
  {\p\over \p t_{j}^{\k}} H\equiv[L_j^{\mu},H]\mod \O_H,
\end{equation}
\end{cor}
\bpf It is easy to show that the ideal of pseudodifference
operators $D$ such that $D\psi=0$ is $\O_H$. From (\ref{dec25}) it
follows that
$$\left(\p_{t_{j}^{\k}} H-[L_j^{\mu},H]\right)\psi=0.$$
Hence, the right and left hand sides of (\ref{nvd}) are equal in
the factor-ring  $\O/\O_H$
\epf

It will be shown below that the system of non-linear equations
(\ref{nvd}) can be regarded as a discrete analog of
Novikov-Veselov hierarchy. The basic equation of this hierarchy --
the discrete analog of the Novikov-Veselov equation
--- is given by (\ref{nvd}) for $j=1$. The operator $L_1^{(1)}$ is
of the form
\begin{equation}\label{nvd1}
  L_1^{(1)}=v(T_1-T_1^{-1}).
\end{equation}
Equation (\ref{nvd}) is equivalent to the system of two equations
for the two functions $u=u_{n,m}(t), \ v=v_{n,m}(t)$:
\begin{equation}\label{nvd2}
  v({\bf t}_1^{-1}u)=u({\bf t}_2v)
\end{equation}
\begin{equation}\label{nvd3}
  \p_t u=\left[({\bf t}_1{\bf t}_2v)({\bf t}_1u)-u({\bf t}_2v)\right]u-[{\bf t}_1{\bf t}_2v-v]
\end{equation}

\smallskip
\noindent {\bf The discrete Novikov-Veselov hierarchy}.

\medskip
The discrete analog of the Novikov-Veselov hierarchy is of
an independent interest. In what follows we consider only the part
of the hierarchy corresponding to ``times'' $t_j:=t_j^1$, and set all
$t_j^2=0$.

Let us write out this part of the hierarchy in a closed
form. We think of it as a system of evolution equations on the following space
\begin{equation}\label{ch1}
  \SP:=\lbrace H,\L\ |\ H=T_1T_2-u(T_1-T_2)-1, \ \L= \sum_{i=0}^{\infty}v_iT_1^{-i+1}\rbrace
\end{equation}
satisfying
\begin{equation}\label{ch2}
  [H,\L]\equiv 0 \mod \O_H^+,
\end{equation}
and such that moreover $u$ and $v_0$ are of the form
\begin{equation}\label{ansatz}
  u=C{({\bf t}_1\tau)\,({\bf t}_2\tau)\over ({\bf t}_1{\bf t}_2 \tau)\, \tau},\ \
  v_0={({\bf t}_1\tau)\,({\bf t}_1^{-1}\tau)\over \tau^2},
\end{equation}
where $C$ is a constant and $\tau=\tau(n,m)$ is some function.

The meaning of (\ref{ch2}) is as follows. A priori the operator
$[H,\L]$ has a unique representation of the form
$$
  [H,\L]=\left(\sum_{s=0}^{\infty}h_sT_1^{-s+2}\right)+DH,
$$
with $D\in \O_1^+$. Therefore, the constraint (\ref{ch2}) is
equivalent to equations $h_s=0$. The first of these equations
$h_0=0$ is an equation for $u$ and $v_0$, which is automatically
satisfied due to (\ref{ansatz}).

By a direct computation of the series expansion of $[H,\L]$ it is
easy to see that equations $h_s=0$ for $s>0$ have the form
\begin{equation}\label{hs0}
  ({\bf t}_2v_s)({\bf t}_1^{-s}u)-({\bf t}_1^{-1}u)v_s=R_s(\tau,v_1,\ldots,v_{s-1}),
\end{equation}
where $R_s$ is some difference polynomial. They recurrently define
$v_s(n,m)$, if ``the initial data'' $v_s|_{m=0}$ are fixed.
Therefore, the space $\SP$ of operators $H,\L$ with the leading
coefficients $u,v_0$ of the form (\ref{ansatz}) satisfying
(\ref{ch2}) can be identified with the space of one function of two
variables and infinite number of functions of one variable, i.e.
$\{\tau(n,m), v_s(n), \, s>0\}$.

Our next goal is to define on $\SP$ a hierarchy of commuting flows.
Any operator in  $\O_{1,0}^+$, and in particular $\L^j$, has a
unique  representation in the form
\begin{equation}\label{ch3}
  \L^j=\sum_{i=-\infty}^{j-1}f_{ij}T_1^i(T_1-T_1^{-1})
\end{equation}
Then the formula (\ref{Li}) with $\k=1$ defines a unique operator
$L_j:=L_j^{(1)}$ such that (\ref{d+}) holds, and also satisfying
the condition (\ref{adj}) with $\k=1$ for the adjoint.
\begin{theo} The equations
\begin{equation}\label{ch9}
  \p_{t_j}\L=[L_j,\L],\quad
  \p_{t_j}H\equiv [L_j,H]\mod\O_H
\end{equation}
define commuting flows on the space $\SP$.
\end{theo}
\bpf Note that the highest power of $T_1$ in $\L$ is $T_1$, and
$\p_{t_j}H=(\p_{t_j}\,u)(T_1-T_2)$. Thus in order to show that
equations (\ref{ch9}) are well-defined we need to prove the
following

(a) $[L_j,\L]$ is of degree not greater than 1;

(b) $[L_j,H]\equiv a_j(T_1-T_2)\mod\O_H$;

(c) the corresponding equations for $v_0$ and $u$ are consistent
with the ansatz (\ref{ansatz}).

\smallskip
The proof of (a) is standard. We compute
\begin{equation}\label{ch5}
  L_j=\L^j+F_j+F_j^1T_1^{-1}+O(T^{-2})\,,
\end{equation}
where
\begin{equation}\label{ch7}
  F_j={\bf t}_1^{-1}f_{1,j}-f_{-1,j},\ F_j^1={\bf t}_1^{-2}f_{2,j}-f_{-2,j}.
\end{equation}
Using $[\L,\L^j]=0$, we get
$$
  [L_j,\L]=[F_j+O(T_1^{-1}),\L]=(F_j-{\bf t}_1F_j)v_0 T_1+O(1),
$$
thus proving (a). Note also that by comparing the leading coefficients we obtain
\begin{equation}\label{ch80}
  {\p\over \p t_j} \ln v_0=F_j-{\bf t}_1F_j
\end{equation}

The proof of (b) is much harder. The difference operator $HL_j$ is of order 1 in $T_2$.
Hence it has a unique representation of the form
\begin{equation}\label{repr}
  HL_j=D_1-a_jT_2+DH,\
\end{equation}
where $D\in \O$ and $D_1\in \O_{1,0}$.

Our next goal is to show that $D_1$ is of degree 1 in $T_1$, i.e. has the form $D_1=b_jT_1+c_j$.
From the equation $T^{-1}_1 H\equiv 0\mod \O_H$ we get
\begin{equation}\label{t2+}
  T_2={\bf t}_1^{-1}u+T_1^{-1}-{\bf t}_1^{-1}u T_1^{-1}T_2=
  {\bf t}_1^{-1}u +(1-{\bf t}_1^{-1}u {\bf t}^{-2}u)T^{-1}+O(T_1^{-2}).
\end{equation}
Equations $[\L^j,H]=0$ and (\ref{ch5}) imply that in $\O^+_H$
the left hand side of (\ref{repr}) is equal to
$$ HL_j=(({\bf t}_1{\bf t}_2F_j-{\bf t}_1F_j)u )T_1+$$
\begin{equation}\label{repr+}
  +\left((1-u{\bf t}_1^{-1}u){\bf t}_1{\bf t}_2F_j+(u {\bf t}^{-1}_1u) {\bf t}_2F_j-F_j+
  ({\bf t}_1^{-1}u){\bf t}_1{\bf t}_2F_j^1-u{\bf t}_1F_j^1\right)
  +O(T^{-1})
\end{equation}
Substituting this expression and the formula for $T_2$ in
(\ref{repr}), we get $D_1=b_jT_1+c+O(T_1^{-1})$, where
\begin{equation}\label{bj}
  b_j:=({\bf t}_1{\bf t}_2F_j-{\bf t}_1F_j)\,u,
\end{equation}
\begin{equation}\label{cj}
  c_j:=a_j{\bf t}_1^{-1}u+(1-u{\bf t}_1^{-1}u){\bf t}_1{\bf t}_2F_j+(u {\bf t}^{-1}_1u)
  {\bf t}_2F_j-F_j+({\bf t}_1^{-1}u){\bf t}_1{\bf t}_2F_j^1-u{\bf t}_1F_j^1
\end{equation}
Now we are going to compute the left and the right hand sides of (\ref{repr}) in
$\O_H^-$. Indeed, in $\O_1^-$ we have
\begin{equation}\label{t2-}
  L_j=-\widetilde\L^j-\widetilde F_j-\widetilde F_j^1T_1+O(T_1^2),
\end{equation}
where, as before,  $\tilde \L=(T_1-T_1^{-1})^{-1} \L^* (T_1-T_1^{-1})$. If $f_{ij}$ are
coefficients of $\wt L$ in (\ref{ch3}), then
\begin{equation}\label{l1}
  \wt \L^j=-\sum_{i=-\infty}^{j-1} T_1^{-i}\cdot
  f_{ij}\cdot (T_1-T_1^{-1}).
\end{equation}
Hence,
\begin{equation}\label{ffs}
  \wt F_j={\bf t}_1F_j=f_{1,j}-{\bf t}_1f_{-1,j},\ \wt
  F_j^1={\bf t}_1^2F_j^1=f_{2,j}-{\bf t}_2f_{-2,j}.
\end{equation}
In order to proceed we now need the following statement.
\begin{lem}
If (\ref{ch2}) is satisfied, then the equation
\begin{equation}\label{ch2dual}
  [H,\wt \L]\equiv 0\mod\O_H^-
\end{equation}
holds.
\end{lem}
\bpf We will prove the lemma by inverting the arguments used above
in the proof of Lemma 3.2. First, for a pair of operators $\L$ and
$H$ satisfying (\ref{ch2}) we introduce a formal solution
$\psi=\psi_{nm}$ of equations
\begin{equation}\label{more1}
  \L\psi=k\psi,\ \ H\psi=0
\end{equation}
of the form
\begin{equation}\label{more2}
  \psi_{nm}=k^{n}\left(\sum_{s=0}^{\infty}\xi_s(n,m)k^{-s}\right).
\end{equation}
Substitution of (\ref{more2}) into (\ref{more1}) gives
a system of difference equations, which recurrently define $\xi_s$.
They have the form
\begin{equation}\label{more3}
  (T_2\,\xi_{s+1})-u\,\xi_{s+1}=\xi_s-u\,(T_2\,\xi_s), \ \
  v_0\,(T_1\,\xi_{s+1})-\xi_{s+1}= \tilde R_s,
\end{equation}
where $\tilde R_s$ are explicit expression linear in the coefficients $v_r$ of $\L$
and difference polynomial in $\xi_r,\ r<s$. If $u,v_0$ are of the form (\ref{ansatz}),
then the first equation for $s=-1$ is satisfied by
\begin{equation}\label{ansatz2}
  \xi_0={{\bf t}_1^{-1}\tau\over \tau}\,.
\end{equation}
The compatibility condition of equations (\ref{more3}) is equivalent to (\ref{ch2}).
These equations uniquely define $\xi_{s+1}$ for all $(n,m)$, if the initial data
$\xi_{s+1}(0,0)$ for (\ref{more3}) is fixed.
Therefore, the solution $\psi$ is unique up to multiplication by a
$(n,m)$-independent Laurent series in the variable $k$.

The function $\psi$ defines a unique operator $\Phi$ of the form (\ref{decphi})
such that equation (\ref{dec313}) holds (with $t_i=0$). Now we define a formal series
\begin{equation}\label{more4}
  \psi^{\s}=k^{-n}\left(\sum_{s=0}^{\infty}\xi^{\s}_s(n,m)
  k^{-s}\right), \ \xi_0^\s={{\bf t}_1\tau\over \tau}
\end{equation}
by the formula
\begin{equation}\label{more5}
  \psi^{\s}=\left((T_1-T_1^{-1})\ \Phi^{-1}
  (T_1-T_1^{-1})^{-1}\right)^*k^{-n}.
\end{equation}
This formal series is an eigenfunction of the operator $\wt L$, i.e. $\wt \L\psi^{\s}=k\psi^{\s}$.
Therefore, in order to prove (\ref{ch2dual}) it is sufficient to prove that $H\psi^{\s}=0$.

From equations (\ref{ansatz},\,\ref{ansatz2}) it follows that
\begin{equation}\label{sdual}
  \wt \psi^{\s}:=H\psi^{\s}=k^{-n}\left(\sum_{s=1}^{\infty} \wt
  \xi^{\s}(n,m)k^{-s}\right)
\end{equation}
Hence, to prove that $\tilde \psi^{\s}=0$ it is enough  to show that
\begin{equation}\label{dec3111}
  \left[\wt \psi^{\s} (T_1^j \psi)\right]_R:=
  \res_k\left({\tilde \psi^{\s}(T_1^j\psi)\,dk\over k^2-1}
  \right)=0, \ \ \forall j\,\geq 2.
\end{equation}
From the definition of $\psi^{\s}$ it follows that
\begin{equation}\label{ind1}
  \left[\psi^{\s} ({\bf t}_1^{2j} \psi)\right]_R=0,\ \
  \left[ \psi^{\s} ({\bf t}_1^{2j+1} \psi)\right]_R=1, \ \ j\geq 0
\end{equation}
(compare to (\ref{dec311}, \ref{dec312})). Using the equation
$H\psi=0$, we get
\begin{equation}\label{ind13}
  {\bf t}_2 \left[\psi^{\s} t^{2j}\psi\right]_R=
  ({\bf t}_1^{2j-1}u)\,\left[{\bf t}_2\psi^{\s}t^{2j-1} \psi)\right]_R-
  ({\bf t}_1^{2j-1}u)\,{\bf t}_2\left[\psi^{\s} t^{2j-1} \psi\right]_R+
  \left[{\bf t}_2\psi^{\s}t^{2j-1} \right]_R
\end{equation}
Then, by induction, it is easy to show that (\ref{ind1}) and (\ref{ind13}) imply
\begin{equation}\label{ind2}
  \left[{\bf t}_2 \psi^{\s}t^{2j+2} \psi)\right]_R=1-
  \prod_{i=0}^{2j+1}({\bf t}_1^iu)^{-1},\quad
  \left[{\bf t}_2 \psi^{\s}t^{2j+1} \psi\right]_R=
  \prod_{i=0}^{2j}({\bf t}_1^iu)^{-1}, \ \ j\geq 0.
\end{equation}
Direct substitution of (\ref{ind2}) into (\ref{dec3111}) completes the proof of the
lemma.
\epf
Now we compute both sides of (\ref{repr}):
\begin{equation}\label{t2++}
  T_2\equiv {1\over u}+\left(1-{1\over u{\bf t}_1u}\right)T^1_1+O(T_1^2) \mod \O_H^-
\end{equation}
Equations (\ref{t2-}) and (\ref{ch2dual})
imply $[L_j,H]=H(\wt F_j+\wt F_j^1 T_1+O(T_1^2))\in \O_H^-$.
Therefore, the operator $D_1$ in (\ref{repr}) has no negative powers of $T_1$,
and hence, it is indeed of the form $b_jT_1+c_j$.

Straightforward computations of the first two coefficients of $[L_j,H]$
give the following formulae
\begin{equation}\label{aj}
  c_j-{a_j\over u}=\wt F_j-{\bf t}_2\wt F_j
\end{equation}
\begin{equation}\label{cj1}
  \left(1-{1\over u{\bf t}_1u}\right)a_j-b_j={1\over {\bf t}_1u}\left({\bf t}_1{\bf t}_2\wt
  F_j+(u{\bf t}_1u-1){\bf t}_2\wt F_j-(u{\bf t}_1u)\,{\bf t}_1\wt F_j-{\bf t}_1u \wt
  F_j^1+u\,{\bf t}_2\wt F_j^1\right)
\end{equation}
From (\ref{bj}, \ref{ffs}) and (\ref{aj}) we get the equations
\begin{equation}\label{abc}
  c_ju=(a_j-b_j)
\end{equation}
and then
\begin{equation}\label{abc1}
  c_j(u{\bf t}_1u-1)={\bf t}_1{\bf t}_2\wt F_j+u{\bf t}_1u\left({\bf t}_2\wt
  F_j-{\bf t}_1\wt F_j\right)-{\bf t}_1u \wt F_j^1+
  u\,{\bf t}_2\wt F_j^1-\wt F_j.
\end{equation}
In order to complete the proof of (b) it is enough now to show that
the right hand side of (\ref{abc1}) is zero. For that we need the
following
\begin{lem} The equations
\begin{equation}\label{J0}
  \wt \F:=-k+(k^2-1)\sum_{j=1}^{\infty}\wt
  F_j\,k^{-j-1}=({\bf t}_1\psi^\s )\psi-
  \psi^\s ({\bf t}_1\psi),
\end{equation}
\begin{equation}\label{J1}
\wt \F^1:=-{({\bf t}_1\tau)^2\over \tau {\bf t}_2\tau}\,k^2+
(k^2-1)\sum_{j=1}^{\infty}\wt F_j^1\,k^{-j-1}=
\psi^\s ({\bf t}_1^2\psi)-({\bf t}_1^2\psi^\s)\psi
\end{equation}
hold.
\end{lem}
\bpf
The expression for the leading coefficients of $\wt \F_j$ and $\wt \F_j^1$ follows from
(\ref{ansatz2}) and (\ref{more5}). In order to prove (\ref{J0}) we need to show
that
\begin{equation}\label{gen1}
  \wt F_j=\res_k\left(\left[({\bf t}_1\psi^\s )\psi-\psi^\s
  ({\bf t}_1\psi)\right]\,{k^jdk\over k^2-1}\right)=
  \res_k\left(\left[(T_1\psi^\s )(\L^j\psi)-\psi^\s
  (T_1\L^j\psi)\right]{dk\over k^2-1}\right)
\end{equation}
From (\ref{more5}), using the relation (\ref{dec315}), we see that
the right hand of (\ref{gen2}) is equal to
\begin{equation}\label{gen2}
  \res_T\left((\L^jT_1^{-1}-T_1\L^j)(T_1-T_1^{-1})^{-1}\right)=f_{1,j}-{\bf t}_1f_{-1,j}
\end{equation}
which proves (\ref{J0}). The proof of (\ref{J1}) is identical.
\epf

From (\ref{J0}) and the equation $H\psi=0$ it follows that
\begin{equation}\label{J3}
  {\bf t}_2\wt \F=-\A\left\{({\bf t}_2\psi^\s(u{\bf t}_1\psi-u{\bf t}_2\psi+\psi)\right\}=
  -\A\left\{({\bf t}_2\psi^\s(u{\bf t}_1\psi+\psi)\right\}.
\end{equation}
Here and below $\A\{\cdot\}$ stands for the antisymmetrization
of the corresponding expression with respect to the interchange of
$\psi^\s$ and $\psi$.

In the same way we get
$$
  {\bf t}_1{\bf t}_2\wt \F=-\A\left\{(u{\bf t}_1\psi^\s-u{\bf t}_2\psi^s+\psi^\s){\bf t}_1^2{\bf t}_2\psi\right\}=$$
$$
-(u{\bf t}_1u){\bf t}_1\A\left\{\psi^\s({\bf t}_1\psi-{\bf t}_2\psi)\right\}-u{\bf t}_2\wt \F^1
-\A\left\{\psi^\s({\bf t}_1(u{\bf t}_1\psi-u{\bf t}_2\psi+\psi)\right\}
$$
Further direct use of the equation $H\psi=0$ and (\ref{J3}) finally gives
the equation
\begin{equation}\label{fin1}
  {\bf t}_1{\bf t}_2\wt \F+u{\bf t}_1u\left({\bf t}_1\wt \F-{\bf t}_2\wt \F\right)-
  {\bf t}_1u \wt \F^1+u\,{\bf t}_2\wt \F^1-\wt\F=0
\end{equation}
The proof of $(b)$ is complete. The comparison of the coefficients  at $T_1$ in the left
and the right hand sides of (\ref{ch9}) gives
\begin{equation}\label{equ}
  \p_{t_j}\ln u =b_j={\bf t}_2 \wt F_j-\wt F_j\end{equation}

Now we are going to prove $(c)$ and derive the evolution equation for $\tau$.
The left and  right action of pseudo-difference operators are formally adjoint,
i.e., for any two operators the equality
$\left(k^{-x}\D_1\right)\left(\D_2k^{x}\right)=
k^{-n}\left(\D_1\D_2k^n\right)+(T_1-1)\left(k^{-x}\left(\D_3k^x\right)\right)$
holds. Here $\D_3$ is a pseudo-difference operator whose coefficients are difference
polynomials in the coefficients of $\D_1$ and $\D_2$. Therefore, from (\ref{J0})
and (\ref{more5}) it follows that
\begin{equation}\label{z8}
\wt \F^{\,0}=-k-(T_1-1)\left((k^2-1)\sum_{s=2}^{\infty}Q_jk^{-j}\right)
\end{equation}
where the coefficients of the series $Q$ are difference polynomials
in the coefficients of the wave operator $\Phi$. Equation (\ref{z8})
implies that
\begin{equation}\label{q81}
\wt F_j=(1-T_1)Q_j=Q_j-{\bf t}_1Q_j.
\end{equation}
Taking into account the ansatz (\ref{ansatz}), we see that equations
(\ref{ch80}) and (\ref{q81}) are equivalent to one equation for the
function $\tau$
\begin{equation}\label{eqntau}
\p_{t_j}\ln \tau = Q_j,
\end{equation}
\begin{rem}
It is necessary to mention that the $Q_j$ are defined only up to an
additive term that is invariant under $T_1$. This ambiguity reflects
the fact that the ansatz (\ref{ansatz}) is invariant under the
transformation
$$\tau(n,m)\longmapsto  f(m) \tau(n,m)$$
where $f(m)$ is an arbitrary function.
\end{rem}
Equation (\ref{eqntau}) completes the proof of the statement that
equations (\ref{ch9}) are well-defined. The proof of the statement
that the corresponding flows on $\SP$ commute with each other is
standard. \epf

\section{Bloch (quasi-periodic) wave solutions.}
To begin with let us prove the implication $(A)\Rightarrow (C)$ in
the main theorem. As it was mentioned above this does not require
the knowledge of the explicit theta-functional form of the function
$\psi$. For the first time an implication of this kind was proved in
\cite{flex}.

Throughout this section $\nu=0,1$ and is considered as an element of
the group $\bbZ_2=\bbZ/2\bbZ$.
\begin{lem}\label{firstlemma}
Let $V\in\bbC^d$, and let $\tau_n^{\nu}(z)$ for
$n\in\bbN,\nu\in\bbZ_2$ be two sequences of holomorphic functions on
$\bbC^d$ such that each divisor $\TC_n^{\nu}:=\lbrace z\in\bbC^d :
\tau_n^{\nu}(z)=0\rbrace $ is not invariant as a set under the shift
by $V$, i.e. $\TC_n^{\nu}\neq \TC_n^{\nu}+V$. Suppose that the
system of equations (considered as a joint system for $\nu=0$ and
$\nu=1$, intertwining $\psi^0$ and $\psi^1$)
\begin{equation}\label{laxdd1}
  \psi_{n+1}^{\nu}(z+V)-u_{n}^{\nu}(z)
  \left(\psi_{n+1}^{\nu+1}(z)-\psi_{n}^{\nu+1}(z+V)\right)-\psi_{n}^{\nu}(z)=0,
\end{equation}
where
\begin{equation}\label{udd1}
  u_n^{\nu}(z)=C\,{\tau_{n+1}^{\nu+1}(z)\,\tau_n^{\nu+1}(z+V)\over
  \tau_{n+1}^{\nu}(z+V)\,\tau_n^{\nu}(z)}
\end{equation}
has solutions $\psi_n^{\nu}$ of the form
\begin{equation}\label{psiddd}
  \psi_n^{\nu}(z)={\a_n^{\nu}(z)\over \tau_n^{\nu}(z)}
\end{equation}
where $\a_n^{\,\nu}$ is a holomorphic function. Then the equation
\begin{eqnarray}\label{taud}
  &\tau_{n+1}^{\nu+1}(z_n^{\nu})\,\tau_n^{\nu+1}(z_n^{\nu}+V)\,\,
  \tau_{n-1}^{\nu}(z_n^{\nu}-V)+
  \tau_{n+1}^{\nu}(z_n^{\nu}+V)\,\tau_n^{\nu+1}(z_n^{\nu}-V)\,\,
  \tau_{n-1}^{\nu+1}(z_n^{\nu})\nonumber\\
  =&\left(\tau_{n+1}^{\nu+1}(z_n^{\nu})\,\tau_n^{\nu+1}(z_n^{\nu}-V)\,\,
  \tau_{n-1}^{\nu}(z_n^{\nu}+V)+
  \tau_{n+1}^{\nu}(z_n^{\nu}-V)\,\tau_n^{\nu+1}(z_n^{\nu}+V)\,\,
  \tau_{n-1}^{\nu+1}(z_n^{\nu})\right)C^{\,2}
\end{eqnarray}
is valid $\forall \ n,\nu,\ \forall z_n^{\nu}\in\TC_n^{\nu}$.
\end{lem}
\bpf Let $I_n^{\nu}(z)$ be the left hand side of (\ref{laxdd1}). A
priori it may have poles at the divisors $\TC_n^{\nu}$ and
$\TC_{n+1}^{\nu}-V$. The vanishing of the residue of $I_n^{\nu}$ at
$\TC_n^{\nu}$ implies
\begin{equation}\label{r1}
  \psi_{n+1}^{\nu+1}(z_n^{\nu})-\psi_n^{\nu+1}(z_n^{\nu}+V)=
  -\a_{n}^{\nu}(z_n^{\nu}){\tau_{n+1}^{\nu}(z_n^{\nu}+V)\over
  \tau_{n+1}^{\nu+1}(z_n^{\nu})\tau_n^{\nu+1}(z_n^{\nu}+V)}\, C^{-1},
\end{equation}
while the vanishing of the residue of $I_{n-1}^{\nu}$ at
$\TC_{n-1}^{\nu}-V$ implies
\begin{equation}
  \psi_{n}^{\nu+1}(z_n^{\nu}-V)-\psi_{n-1}^{\nu+1}(z_n^{\nu})=
  \a_{n}^{\nu}(z_n^{\nu}){\tau_{n-1}^{\nu}(z_n^{\nu}-V)\over \tau_{n}^{\nu+1}(z_n^{\nu}-V)
  \tau_{n-1}^{\nu+1}(z_n^{\nu})}\,  C^{-1}. \label{r2}
\end{equation}
On the other hand, the evaluation of $I_n^{\nu+1}$ at the divisor
$\TC_n^{\nu}-V$ implies
\begin{equation}\label{r3}
  \psi_{n+1}^{\nu+1}(z_n^{\nu})-\psi_{n}^{\nu+1}(z_n^{\nu}-V)=
  -\a_{n}^{\nu}(z_n^{\nu}){\tau_{n+1}^{\nu}(z_n^{\nu}-V)\over \tau_{n+1}^{\nu+1}(z_n^{\nu})
  \tau_{n}^{\nu+1}(z_n^{\nu}-V)}\, C,
\end{equation}
while the evaluation of $I_{n-1}^{\nu+1}$ at the divisor
$\TC_n^{\nu}$ implies
\begin{equation}\label{r4}
  \psi_{n}^{\nu+1}(z_n^{\nu}+V)-\psi_{n-1}^{\nu+1}(z_n^{\nu})=
  \a_{n}^{\nu}(z_n^{\nu}){\tau_{n-1}^{\nu}(z_n^{\nu}+V)\over \tau_{n}^{\nu+1}(z_n^{\nu}+V)
  \tau_{n-1}^{\nu+1}(z_n^{\nu})}\,  C.
\end{equation}
The left-hand-side of the difference of (\ref{r1}) and
(\ref{r2}) is the same as that of the difference of (\ref{r3}) of
(\ref{r4}); equating the right-hand-sides of these differences
yields (\ref{taud}). \epf

Formulation $(A)$ of our main theorem implies that the assumption
of lemma \ref{firstlemma} is satisfied for $C=c_3$, $z\in \bbC^g$,
and
\begin{equation}\label{tauexplicit}
  \tau_n^\nu(z)=\theta\left(Un+(1-\nu)W+z\right)\,
  \left(c_1^{(\,l,\,z)}c_2^n\right)^{\nu-\frac12},
\end{equation}
where $l\in\bbC^g$ is a vector such that $(l,V)=1$. Then from
(\ref{taud}) for $\nu=0$ we get on the divisor $\TC_0^0$, i.e. for
$\theta(Z)=\theta(z+W)=0$,
\begin{eqnarray}\label{95}
  \tau_1^1(z)\,\tau_0^1(z+V)\,\tau_{-1}^0(z-V)+\tau_1^0(z+V)\,
  \tau_0^1(z-V)\,\tau_{-1}^1(z)\nonumber\\
  =c_3^2(\tau_1^1(z)\,\tau_0^1(z-V)\,\tau_{-1}^0(z+V)+
  \tau_1^0(z-V)\,\tau_0^1(z+V)\,\tau_{-1}^0(z))\,,
\end{eqnarray}
which upon substituting (\ref{tauexplicit}) yields, after canceling
the common factors,
\begin{eqnarray}\nonumber
  c_1^2c_2^2\,\theta(Z+U-W)\,\theta(Z+V-W)\,\theta(Z-U-V)&\\
  +\theta(Z+U+V)\,\theta(Z-V-W)\,\theta(Z-U-W)&\nonumber\\
  =c_2^2c_3^2\,\theta(Z+U-W)\,\theta(Z-V-W)\,\theta(Z+V-U)&\nonumber\\
  +c_1^2c_3^2\,\theta(Z-V+U)\,\theta(Z+V-W)\,\theta(Z-U-W)
\end{eqnarray}
which is identical to equation (\ref{cm7d}) with the minus sign
chosen for $W$ (and correspondingly the constants $c_1$ and $c_2$
appearing in positive power). Similarly the case of $\nu=1,n=0$ of
formula (\ref{taud}) yield the plus sign case of (\ref{cm7d}). The
implication $(A)\Rightarrow (C)$ in the main theorem is thus proved.

\bigskip
Let us now show that $(C)$ can also be obtained as a corollary of a
more general {\it fourth order} relation for Prym theta-functions.
As it was mentioned above, in \cite{grin} it was proved that
equation (\ref{laxdd}) implies the five-term equation (\ref{51}).
Note, that all the pairs of indices have sums of the same
parity, i.e. equation (\ref{51}) is in fact a pair of equations on
two functions $\psi$ defined on two sublattices of the variables
$(n,m)$.

The statement that $\psi_{n,m}$ satisfy (\ref{51}) can be proved
directly. Indeed all the functions involved in the equation are in
$$
  H^0(D+(n+1)P_1^+-(n-1)P_1^-+(m+1)P_2^+-(m-1)P_2^-+\nu(P_3^+-P_3^-))
$$
By the Riemann-Roch theorem the dimension of the latter space is $4$. Hence,
any five elements of this space are linearly dependent, and it remains to
find the coefficients of (\ref{51}) by a
comparison of singular terms at the points $P_1^{\pm}, P_2^{\pm}$.
For $n+m=0\mod 2)$ we get
\begin{eqnarray}\label{52}
  \tilde a_{n,m}&=&c_1^2c_3^2\,{\theta(Z_{n,m}+V)\,\theta(Z_{n,m}+U-V+W)\over
  \theta(Z_{n,m}-V)\,\theta(Z_{n,m}+U+V+W)},\nonumber\\
  \tilde b_{n,m}&=&c_2^2c_3^2\,{\theta(Z_{n,m}+U)\,\theta(Z_{n,m}-U+V+W)\over
  \theta(Z_{n,m}-U)\,\theta(Z_{n,m}+U+V+W)},\nonumber\\
  \tilde c_{nm}&=&c_1^2c_2^2\,{\theta(Z_{n,m}+U)\,\theta(Z_{n,m}+V)\,\theta(Z_{n,m}-U-V+W)\over
  \theta(Z_{n,m}+U)\,\theta(Z_{n,m}-V)\,\theta(Z_{n,m}+U+V+W)},
\end{eqnarray}
where $Z_{n,m}=Z+Un+Vm$. From the normalization of $\psi_{n,m}$ it
follows that
\begin{equation}\label{54}
  \tilde d_{nm}=1-\tilde a_{n,m}-\tilde b_{n,m}+\tilde c_{n,m}
\end{equation}
Substituting here (\ref{psipr}) and (\ref{52}-\ref{54})
proves the following statement.
\begin{lem}
For any four points $A,U,V,W$ on the image
$\G\hookrightarrow\P(\G)$, and any $Z\in\P(\G)$ the following
equation holds: \hfill

\begin{eqnarray}\label{quad}
  \theta(Z+W)\times [ &\theta(A+U+V+Z)\,\theta(Z-U)\,\theta(Z-V)\nonumber\\
 &-c_1^2c_3^2\,\theta(A+U-V+Z)\,\theta(Z-U)\,\theta(Z+V)\nonumber\\
  &-c_2^2c_3^2\,\theta(A-U+V+Z)\,\theta(Z+U)\,\theta(Z-V)\nonumber\\
 &+c_1^2c_2^2\,\theta(A-U-V+Z)\,\theta(Z+U)\,\theta(Z+V)]= \nonumber \\
  =\theta(A+Z)\times [ &\theta(W+U+V+Z)\,\theta(Z-U)\,\theta(Z-V)\nonumber\\
 &-c_1^2c_3^2\,\theta(W+U-V+Z)\,\theta(Z-U)\,\theta(Z+V)\nonumber\\
 &-c_2^2c_3^2\,\theta(W-U+V+Z)\,\theta(Z+U)\,\theta(Z-V)\nonumber\\
 &+c_1^2c_2^2\,\theta(W-U-V+Z)\,\theta(Z+U)\,\theta(Z+V)].\nonumber\\
\end{eqnarray}
\end{lem}
To the best of the authors' knowledge equation (\ref{quad}) is a new
identity for Prym theta-functions. For $Z$ such that $\theta(W+Z)=0$
it is equivalent to equation (\ref{cm7d}) with minus sign chosen.
The second equation of the pair (\ref{cm7d}) can be obtained from
(\ref{51}) considered for the odd case, i.e. for $n+m=1\, \mod 2$.

\bigskip
\noindent {\bf Wave solutions.} In section 2 we proved that if
$\theta(Z)$ is the Prym theta function, then equation (\ref{laxdd})
with $u$ as in (\ref{ud}) has not just one solution $\psi$ of the
form (\ref{pd}) but a family of them parameterized by points $A$ in
the image $\G\longmapsto \P(\G)$ under the Abel-Prym map. Note,
however, that formulation $(C)$ of the main theorem does not involve
$A$. The first step in proving the ``only if'' part of $(C)$ (and
thus also of (A) and (B), which imply (C)) is to introduce a {\it
spectral parameter} in the problem, i.e. to show that equations
(\ref{cm7d}) are sufficient for the existence of certain formal
solutions of equations (\ref{laxdd1}). These solutions are functions
of the form
\begin{equation}\label{w1}
  \psi_n^{\nu}(z)=k^n C^{(l,z)}\,\phi_n^\nu(z,k),
\end{equation}
where $k^{-1}$ is a formal parameter (eventually to be identified
with the local coordinate on the curve), $\phi_n^\nu(z,k)$ is a
regular series in $k^{-1}$, i.e.
\begin{equation}\label{w2}
  \phi_n^{\nu}(z,k)=\sum_{s=0}^{\infty}\xi_{n,\,s}^{\nu}(z)k^{-s},
\end{equation}
and $l\in\bbC^d$ is such that $(l,V)=1$.

The ultimate goal of this section is to show that such solutions
exist with $\xi_{n,s}^\nu$ being holomorphic functions of $z\in
\bbC^g$, defined outside the divisor $\theta (z+Un+(1-\nu)W)=0$
\footnote{In \cite{kr-schot,kr-prym,kr-tri} the
corresponding solutions were called $\l$-periodic reflecting the normalization
leading to their definition. The idea of that normalization goes back to \cite{kp1}.}.
As we shall see below, an obstruction for the existence of such solutions
is the ``bad locus''
$$
 \Sigma:=\Sigma^{\,0}\cup \Sigma^1,
$$
where $\Sigma^\nu$ is the $V$-invariant subvariety of the divisor
$\Theta+(\nu-1)W$ that is not $U$-invariant, i.e.
\begin{equation}\label{def.sigma}
  \Sigma^\nu:=\left\lbrace Z\in X:\quad\begin{array}{l} \forall n\in\bbZ\quad \theta(Z+nV+(1-\nu)W)=0;\\
  \exists n\in\bbZ\quad \theta(Z+U+nV+(1-\nu)W)\ne 0\end{array}\right\rbrace
\end{equation}
We will prove in lemma \ref{no.sigma} that the bad locus is empty,
but until then we construct the wave solutions with the desired
properties only along certain affine subspaces of $\bbC^g$; then
we will patch these together. \footnote{The locus $\Sigma$ is an
analog of singular locus considered in \cite{shiota}. The authors
are grateful to Enrico Arbarello for an explanation of  its
crucial role, which helped them to focus on the heart of the
problem.}

\smallskip
{\bf Notations.} Denote $\pi: \bbC^g\to X=\bbC^g/\Lambda$ the
universal cover map for $X$. Let $Y$ be the Zariski closure of the
group $\langle\bbZ V\rangle\subset X$. As an abelian subvariety, it
is generated by its irreducible component $Y^0$, containing $0$, and
by the point $V_0$ of finite order in $X$, such that $V-V_0\in
Y^0,\, NV_0=\l_0\in \Lambda$. Shifting $Y$ if needed, we may assume,
without loss of generality, that $0$ is not in the bad locus
$\Sigma$. Since any subset of $Y$ that is invariant under the shift
by $V$ is dense in $Y$, this implies that $Y\cap \Sigma=\emptyset$.

We denote $\C:=\pi^{-1}(Y)$. Then $\C$ is a union of its connected
component passing through zero (which is a linear subspace $\mathbb
V\cong\bbC^d\subset\bbC^g$) and shifts by a preimage of a vector of
finite order, i.e. we have $\C=\cup_{r\in \bbZ} (\mathbb V+rV_0)$.
Denoting then $\Lambda_0:=\Lambda\cap {\C}$ we have
$Y=\C/\Lambda_0$, and we can also write $\Lambda_0=\wt
\Lambda_0+\bbZ V_0$, where $\wt \Lambda_0:=\Lambda\cap \mathbb V$.

\smallskip
In what follows we assume that $\tau_n^\nu(z)$ are non-vanishing
identically holomorphic functions of the variable $z\in \C$ having
the following factors of automorphy with respect to $\Lambda_0$:
\begin{equation}\label{mon1}
  \tau_n^{\nu}(z+\l)=\tau_n^{\nu}(z)\,e^{(z,\,\a_\l)+n \b_\l+w_\l^\nu}
\end{equation}
where $\a_\l,\b_\l^{\nu}$ are independent of $n$,
and we define for further use
\begin{equation}\label{blambda}
  b_\l^\nu:=e^{\b_{\l}+w_{\l}^\nu-w_\l^{\nu+1}}.
\end{equation}
This means that $u_n^\nu(z)$ given by (\ref{udd1}) is a section of
some degree zero line bundle on $Y$.

\begin{prop}\label{solutions}
Suppose equation (\ref{taud}) for $\tau_n^\nu(z)$ holds.
Then equations (\ref{laxdd1}) with potentials $u_n^\nu(z)$ given by (\ref{udd1})
have wave solutions of the form (\ref{w1}) such that
\begin{itemize}
\item[(i)] the coefficients $\xi_{n,\,s}^\nu(z)$ of the formal
series $\phi_n^\nu(z,k)$ are meromorphic functions of the
variable $z\in \C$ with a simple pole at the divisor
$\TC_n^\nu$,
\begin{equation}\label{v1}
  \xi_{n,\,s}^\nu(z)={\tau_{n,s}^{\nu+1}(z)\over \tau_n^\nu(z)}\,,
\end{equation}
where $\tau_{n,s}^{\nu+1}(z)$ is a holomorphic function (the
shift from $\nu$ to $\nu+1$ is only for notational ease to simplify
further formulas), and
\begin{equation}\label{v0}
  \tau_{n,0}^{\nu}(z)=\tau_{n-1}^\nu(z).
\end{equation}
\item[(ii)] Each of the individual terms in the power series
expansion of $\phi$ have the following automorphy properties
(note we are not yet making any claims regarding $\phi$ as a
whole)
\begin{equation}\label{bloch2}
 b_\l^\nu \xi_{n,s}^\nu (z+\l)-\xi_{n,s}^{\nu}(z)=
\sum_{i=1}^s B_{i,n-s+i}^\l\xi_{n,s-i}^{\nu}(z),
\end{equation}
for any $\l\in\Lambda_0$ (notice that the coefficients depend on
$i$ and in a sort of diagonal way on $n$, but do {\it not}
depend on $\nu$, which will be important for the future
computations).
\end{itemize}
\end{prop}
\bpf Writing down the equation for $\psi$ in terms of the power
series expansions in $k^{-s}$, and equating coefficient of $k^{-s}$
to zero (i.e. substituting (\ref{udd1},\ref{w1},\ref{w2}) into
(\ref{laxdd1})) yields
\begin{equation}\label{xi.all.s}
  C\,\xi_{n+1,\, s+1}^{\nu}(z+V)-u_{n}^\nu(z)\,(\xi_{n+1,\,s+1}^{\nu+1}(z)-
  C\,\xi_{n,\,s}^{\nu+1}(z+V))+\xi_{n,s}^{\nu}(z)=0.
\end{equation}
For $s=-1$ equation (\ref{xi.all.s}) is satisfied with
$\tau_{n,0}^\nu$ given by (\ref{v0}), i.e. with
\begin{equation}\label{w4}
\xi^{\nu}_{n,\,0}(z)={\tau_{n-1}^{\nu+1}(z)\over \tau_{n}^{\nu}(z)}\,.
\end{equation}
We will now prove the lemma by induction in $s$. Let us assume
inductively that for $r\le s-1$ the functions $\xi_{n,r}^\nu(z)$ are
known, {\it for all $n$ and $\nu$}, and satisfy the
quasi-periodicity condition (\ref{bloch2}) above --- it is customary in the subject to call such solutions Bloch solutions or Bloch functions.

The idea of the proof of the inductive step is as follows. We write
down the equation relating $\tau_{n+1,s+1}^\nu$ (we are using $n+1$
instead of $n$ solely for the ease of notations --- recall that the
inductive assumption is for all $n$) to the $\tau$ for smaller
values of $s$ (which we know inductively to exist and be
holomorphic). From this equation we then get an explicit formula for
$\tau_{n+1,s+1}^\nu$ on the divisor $\TC_n^\nu$, i.e. for
$\tau_n^\nu(z)=0$. We also get an explicit formula for
$\tau_{n+1,s+1}^\nu$ for $z$ such that $\tau_n^\nu(z+V)=0$, which
after translating the argument gives another formula for
$\tau_{n+1,s+1}^\nu$ on the divisor $\TC_n^\nu$. Once we verify that
the two resulting formulas agree (this is a hard computation using
the step of the induction), it will follow that $\tau_{n+1,s+1}^\nu$
restricted to $\TC_n^\nu$ is in fact holomorphic and thus can be
extended from this divisor holomorphically to $\bbC^d$. We now give
the details of this argument.

\smallskip
Writing down equation (\ref{xi.all.s}) in terms of $\tau$'s for
arbitrary $s$, and clearing denominators yields
\begin{eqnarray}\label{tau.all.s}
  C\tau_{n+1,s+1}^{\nu+1}(z+V)\tau_n^\nu(z)
  &-&C\tau_n^{\nu+1}(z+V)\tau_{n+1,s+1}^\nu(z)\nonumber\\
  -C^2\tau_{n,s}^\nu(z+V)\tau_{n+1}^{\nu+1}(z)
  &+&\tau_{n+1}^\nu(z+V)\tau_{n,s}^{\nu+1}(z)=0
\end{eqnarray}
These equations can be easily solved on the divisor $\TC_n^\nu$.
Indeed, if we take $z=z_n^\nu\in\TC_n^\nu$ here, the first term will
vanish, and we get the following formula
\begin{equation}\label{expr1}
  C\tau_{n+1,s+1}^\nu(z_n^\nu)=\frac{
  \tau_{n,s}^{\nu+1}(z_n^\nu)\tau_{n+1}^\nu(z_n^\nu+V)-
  C^2\tau_{n,s}^\nu(z_n^\nu+V)\tau_{n+1}^{\nu+1}(z_n^\nu)}
  {\tau_n^{\nu+1}(z_n^\nu+V)}.
\end{equation}
Alternatively, using equation (\ref{tau.all.s}) for $\nu+1$ instead
of $\nu$ and setting $z=z_n^\nu-V$, for $z_n^\nu\in\TC_n^\nu$ as
above, we get
\begin{equation}\label{expr2}
  C\tau_{n+1,s+1}^\nu(z_n^\nu)=\frac{
  \tau_{n,s}^{\nu+1}(z_n^\nu)
  \tau_{n+1}^\nu(z_n^\nu-V)-C^2\tau_{n,s}^\nu(z_n^\nu-V)
  \tau_{n+1}^{\nu+1}(z_n^\nu)}{\tau_n^{\nu+1}(z_n^\nu-V)}.
\end{equation}
For $\tau_{n+1,s+1}^\nu$ to have a chance to exist, these two
expressions have to agree.
\begin{lem}\label{induct.ok}
If the inductive assumption (and the conditions of the proposition,
in particular formula (\ref{taud})) is satisfied for $s$, then the
two expressions above for the function $\tau_{n+1,s+1}^\nu(z)$
restricted to the divisor $\TC_n^\nu$ are equal.
\end{lem}
\bpf
Equating the two expressions
obtained for $\tau_{n+1,s+1}^\nu$ on $\TC_n^\nu$, we see that what
we need to prove is the following identity
\begin{eqnarray}
  \tau_{n,s}^{\nu+1}(z_n^\nu)\tau_{n+1}^\nu(z_n^\nu-V)\tau_n^{\nu+1}(z_n^\nu+V)
  -C^2\tau_{n,s}^\nu(z_n^\nu-V)\tau_{n+1}^{\nu+1}(z_n^\nu)\tau_n^{\nu+1}(z_n^\nu+V)
  \nonumber\\
  =\tau_{n,s}^{\nu+1}(z_n^\nu)\tau_{n+1}^\nu(z_n^\nu+V)\tau_n^{\nu+1}(z_n^\nu-V)-
  C^2\tau_{n,s}^\nu(z_n^\nu+V)\tau_{n+1}^{\nu+1}(z_n^\nu)\tau_n^{\nu+1}(z_n^\nu-V).
  \label{needed}
\end{eqnarray}
To prove that this is the case, we will use the inductive assumption
for $n-1,s-1$, and equation (\ref{taud}). Indeed, for $n-1,s-1$
equation (\ref{tau.all.s}) reads
\begin{eqnarray}
  C\tau_{n,s}^{\nu+1}(z+V)\tau_{n-1}^\nu(z)
  &-C\tau_{n-1}^{\nu+1}(z+V)\tau_{n,s}^\nu(z)&\nonumber\\
  -C^2\tau_{n-1,s-1}^\nu(z+V)\tau_n^{\nu+1}(z)
  &+\tau_n^\nu(z+V)\tau_{n-1,s-1}^{\nu+1}(z)&=0\nonumber
\end{eqnarray}
By the inductive assumption we know that this is satisfied. If we
now take $z=z_n^\nu-V$, i.e. set $\tau_n^\nu(z+V)=0$ here, we
get
$$
  C^2\tau_{n-1,s-1}^\nu(z_n^\nu)=\frac{
  C\tau_{n,s}^{\nu+1}(z_n^\nu)\tau_{n-1}^\nu(z_n^\nu-V)
  -C\tau_{n-1}^{\nu+1}(z_n^\nu)\tau_{n,s}^\nu(z_n^\nu-V)}
  {\tau_n^{\nu+1}(z_n^\nu-V)}.
$$
Similarly, if we instead take the equation with $\nu+1$ instead of
$\nu$, and take $z=z_n^\nu$, we get
$$
  \tau_{n-1,s-1}^\nu(z_n^\nu)=\frac{
  C\tau_{n-1}^\nu(z_n^\nu+V)\tau_{n,s}^{\nu+1}(z_n^\nu)
  -C\tau_{n,s}^\nu(z_n^\nu+V)\tau_{n-1}^{\nu+1}(z_n^\nu)}
  {\tau_n^{\nu+1}(z_n^\nu+V)}
$$
Since we inductively assumed the existence and uniqueness of
$\tau_{n-1,s-1}^\nu$, these two expressions must agree, which is to
say that we have the following identity
\begin{eqnarray}
  \tau_{n,s}^{\nu+1}(z_n^\nu)\tau_{n-1}^\nu(z_n^\nu-V)
  \tau_n^{\nu+1}(z_n^\nu+V)-\tau_{n,s}^\nu(z_n^\nu-V)
  \tau_{n-1}^{\nu+1}(z_n^\nu)\tau_n^{\nu+1}(z_n^\nu+V)\nonumber\\
  =C^2\tau_{n,s}^{\nu+1}(z_n^\nu)\tau_{n-1}^\nu(z_n^\nu+V)
  \tau_n^{\nu+1}(z_n^\nu-V)-C^2\tau_{n,s}^\nu(z_n^\nu+V)
  \tau_{n-1}^{\nu+1}(z_n^\nu)\tau_n^{\nu+1}(z_n^\nu-V)\label{known}
\end{eqnarray}

Notice now how similar this known identity is to formula (\ref{needed})
that we need to prove. Indeed, the coefficient of
$\tau_{n,s}^{\nu+1}(z_n^\nu)$ in (\ref{needed}) is equal to
$$
  \tau_{n+1}^\nu(z_n^\nu-V)\tau_n^{\nu+1}(z_n^\nu+V)
  -\tau_{n+1}^\nu(z_n^\nu+V)\tau_n^{\nu+1}(z_n^\nu-V).
$$
Now using formula (\ref{taud}), which we know holds for $\tau$, we
see that this coefficient is equal to
\begin{equation}
  \left(\tau_n^{\nu+1}(z_n^\nu-V)\tau_{n-1}^\nu(z_n^\nu+V)
  -\tau_n^{\nu+1}(z_n^\nu+V)\tau_{n-1}^\nu(z_n^\nu-V)\right)
  \frac{\tau_{n+1}^{\nu+1}(z_n^\nu)}{\tau_{n-1}^{\nu+1}(z_n^\nu)}
\end{equation}
Substituting this expression into (\ref{needed}) is
equivalent to the identity
\begin{eqnarray}
  \tau_{n,s}^{\nu+1}(z_n^\nu)
  \frac{\tau_{n+1}^{\nu+1}(z_n^\nu)}{\tau_{n-1}^{\nu+1}(z_n^\nu)}
  \left(\tau_n^{\nu+1}(z_n^\nu-V)\tau_{n-1}^\nu(z_n^\nu+V)
  -C^2\tau_n^{\nu+1}(z_n^\nu+V)\tau_{n-1}^\nu(z_n^\nu-V)\right)\\
  =\tau_{n,s}^\nu(z_n^\nu-V)\tau_{n+1}^{\nu+1}(z_n^\nu)\tau_n^{\nu+1}(z_n^\nu+V)\nonumber
  -C^2\tau_{n,s}^\nu(z_n^\nu+V)\tau_{n+1}^{\nu+1}(z_n^\nu)\tau_n^{\nu+1}(z_n^\nu-V).
\end{eqnarray}
Multiplying this identity by
$\frac{\tau_{n-1}^{\nu+1}(z_n^\nu)}{\tau_{n+1}^{\nu+1}(z_n^\nu)}$
yields formula (\ref{known}), which we inductively know to hold.
Thus formula (\ref{needed}) holds, and the lemma is proven. \epf

\begin{lem}
The function $\tau_{n+1,s+1}^\nu(z_n^\nu)$ given by (\ref{expr1})
and (\ref{expr2}) can be extended to a holomorphic function on the
entire divisor $\TC_n^\nu$.
\end{lem}
\bpf
The expression (\ref{expr1}) for
$\tau_{n+1,s+1}^\nu(z_n^\nu)$ is certainly holomorphic when
$\tau_n^{\nu+1}(z_n^\nu+V)$ is non-zero, i.e. is holomorphic
outside of $\TC_n^\nu\cap(\TC_n^{\nu+1}-V)$. Similarly the
expression for $\tau_{n+1,s+1}^\nu$ given by formula
(\ref{expr2}) is holomorphic away
from  $\TC_n^\nu\cap(\TC_n^{\nu+1}+V)$.

We have assumed that the closure of the abelian subgroup generated
by $V$ is everywhere dense. Thus for any $z_n^\nu\in\TC_n^\nu$
there must exist some $N\in\bbN$ such that
$z_n^\nu+(N+1)V\not\in\TC_n^{\nu+1}$; let $N$ moreover be the
minimal such $N$. From (\ref{expr1}) it then follows that
$\tau_{n+1,s+1}^\nu$ can be extended holomorphically to the point
$z_n^\nu+NV$. However, by lemma \ref{induct.ok} we know that the
expressions (\ref{expr2}) and (\ref{expr1}) agree. Thus expression
(\ref{expr2}) must also be holomorphic at $z_n^\nu+NV$; since its
denominator there vanishes, it means that the numerator must also
vanish, i.e. we must have
$$
C\tau_{n,s}^{\nu+1}(z_n^\nu+NV)
  \tau_{n+1}^\nu(z_n^\nu+(N-1)V)-\tau_{n,s}^\nu(z_n^\nu+(N-1)V)
  \tau_{n+1}^{\nu+1}(z_n^\nu+NV)=0.
$$
But this expression is equal to the numerator of (\ref{expr1}) at
$z_n^\nu+(N-1)V$; thus $\tau_{n+1,s+1}^\nu$ defined from
(\ref{expr1}) is also holomorphic at $z_n^\nu+(N-1)V$ (the
numerator vanishes, and the vanishing order of the denominator is
one, since we are talking exactly about points on its vanishing
divisor). Thus unless $N=0$ we have a contradiction, since $N$ was
chosen minimal. For $N=0$, however,
$z_n^\nu+V\not\in\TC_n^{\nu+1}$, and thus (\ref{expr1}) defines
$\tau_{n+1,s+1}^\nu$ holomorphically at $z_n^\nu$.\epf

Recall now that an analytic function on an analytic divisor in
$\bbC^d$ has a holomorphic extension to all of $\bbC^d$
(\cite{serr}). Therefore, there exists a holomorphic function $\wt
\tau_{n+1,s+1}^{\nu}(z)$ extending the function given on the divisor
$\TC_n^\nu$ by the r.h.s. of (\ref{expr1}) (by the above lemma, it
is holomorphic, and thus the extension is holomorphic). It is then
natural to attempt to use the function $\wt
\xi_{n+1,s+1}^{\nu}:=\wt\tau_{n+1,s+1}^{\nu+1}/\tau_{n+1}^\nu$ for
the proposition, but this cannot be done immediately, as such an
extension does not need to be quasi-periodic, nor is going to be a
solution of equation (\ref{xi.all.s}). We will thus need to adjust
this extension appropriately.

We start by determining the quasi-periodicity properties: indeed,
for $z_n^{\nu+1}\in\TC_n^{\nu+1}$, where we know that
$\wt\tau_{n+1,s+1}^{\nu+1}$ is given by (\ref{expr1}), we have
\begin{equation}\label{n.to.n+1}
  \wt\xi_{n+1,s+1}^\nu(z_n^{\nu+1})=-C\xi_{n,s}^\nu
  (z_n^{\nu+1}+V)+\frac{\tau_{n,s}^\nu(z_n^{\nu+1})
  \tau_{n+1}^{\nu+1}(z_n^{\nu+1}+V)}
  {\tau_n^\nu(z_n^{\nu+1}+V)\tau_{n+1}^\nu(z_n^{\nu+1})},
\end{equation}
from which by using the quasi-periodicity of $\tau_n$ (\ref{mon1})
and that of $\tau_{n,s}$ (\ref{bloch2}), it follows that
\begin{eqnarray}
  b_\l^\nu\wt\xi_{n+1,s+1}^\nu(z_n^{\nu+1}+\l)=
  -C\left(\xi_{n,s}^\nu(z_n^{\nu+1}+V)
  -\sum_{i=1}^s B_{i,n-s+i}^\nu\xi_{n,s-i}^\nu(z_n^{\nu+1}+V)\right)
  \nonumber\\
  +\frac{\left(\tau_{n,s}^\nu(z_n^{\nu+1})
  +\sum\limits_{i=1}^s B_{i,n-s+i}^\l\tau_{n,s-i}^\nu(z_n^{\nu+1})\right)
  \tau_{n+1}^{\nu+1}(z_n^{\nu+1}+V)}{\tau_n^\nu(z_n^{\nu+1}+V)
  \tau_{n+1}^\nu(z_n^{\nu+1})}
\end{eqnarray}
since the $e^{(2z+V,\alpha_\l)+(2n+1)\b_\l}$  factors for the
second term of (\ref{n.to.n+1}) coming from (\ref{mon1}) cancel in
the numerator and denominator, and the remaining
$e^{2\omega_\l^{\nu+1}-2\omega_\l^\nu}$ cancels  with
$b_\l^\nu/b_\l^{\nu+1}$. We now note that the terms in the r.h.s.
split in pairs similar to those in (\ref{n.to.n+1}) and we can thus
simplify this to get
\begin{equation}
  0=b_\l^\nu\wt\xi_{n+1,s+1}^\nu(z_n^{\nu+1}+\l)-
  \wt\xi_{n+1,s+1}^\nu(z_n^{\nu+1})-
  \sum\limits_{i=1}^s B_{i,n-s+i}^\l\xi_{n+1,s+1-i}^\nu(z_n^{\nu+1})
\end{equation}
This says that the function on the right-hand-side here --- denote
it by $g_{n+1,s+1}^{\lambda,\nu}(z)$ --- vanishes for
$z=z_n^{\nu+1}\in\TC_n^{\nu+1}$ and has a pole for
$z\in\TC_{n+1}^\nu$. Using formula (\ref{w4}) for $\xi_{n,0}^\nu$,
we can then write
$$
 g_{n+1,s+1}^{\l,\nu}(z)=f_{n+1,s+1}^{\l,\nu}(z)\xi_{n+1,0}^{\nu}(z),
$$
where $f_{n+1,s+1}^{\l,\nu}(z)$ is now holomorphic, and satisfies
the twisted homomorphism relations
\begin{equation}\label{bl3}
f_{n+1,s+1}^{\l+\mu,\nu}(z)=f_{n+1,s+1}^{\l,\nu}(z+\mu)+f_{n+1,s+1}^{\mu,\nu}(z)
\end{equation}
We only know the function $\wt\xi$ to have the desired quasi-periodicity on
the divisor $\TC_n^{\nu+1}$, and would now like to adjust it so that
the corrected function would have computable quasi-periodicity for all $z$.
To achieve this, we need to add to $\wt\xi$ a summand involving $f$.

Indeed, $f$ defines an element of the first cohomology group of
$\Lambda_0$ with coefficients in the sheaf of holomorphic functions,
$f\in H^1_{gr}(\Lambda_0,H^0(\bbC^d, \O))$. The arguments
identical to that in the proof of part (b) of Lemma 12 in
\cite{shiota} show that there must then exist a holomorphic function
$h_{n+1,s+1}^\nu(z)$ such that
\begin{equation}\label{bl4}
  f_{n+1,s+1}^{\l,\nu}(z)=h_{n+1,s+1}^\nu(z+\l)-h_{n+1,s+1}^\nu(z)
  +E_{n+1,s+1}^{\l,\,\nu},
\end{equation}
where $E_{n+1,s+1}^{\l,\,\nu}$ is a ($z$-independent!) constant. By
using equation (\ref{bl3}) we observe that $E$ depends on $\l$
linearly, i.e. that
\begin{equation}\label{bl5}
 E_{n+1,s+1}^{\l+\mu,\nu}=E_{n+1,s+1}^{\l,\nu}+E_{n+1,s+1}^{\mu,\nu}
\end{equation}
We then define
$$
  \zeta_{n+1,s+1}^{\,\nu}(z):=\wt\xi_{n+1,s+1}^{\nu}(z)
  -h_{n+1,s+1}^\nu(z)\xi_{n+1,0}^\nu(z).
$$
Using (\ref{w4}) and (\ref{mon1}), we first compute
\begin{equation}\label{xi.0.mon}
 \xi_{n+1,0}^\nu(z+\l)=\frac{\tau_n^{\nu+1}(z+\l)}{\tau_{n+1}^\nu
 (z+\l)}=e^{\omega_\l^{\nu+1}-\omega_\l^\nu-\b_\l}\frac{\tau_n^{
 \nu+1}(z)}{\tau_{n+1}^\nu(z)}=\frac{\xi_{n+1,0}^\nu(z)}{b_\l^\nu}
\end{equation}
and then compute the quasi-periodicity
\begin{eqnarray}\label{bl4a}
 &b_\l^\nu\zeta_{n+1,s+1}^{\,\nu}(z+\l)-\zeta_{n+1,s+1}^\nu(z)=
 (b_\l^\nu\wt\xi_{n+1,s+1}^{\,\nu}(z+\l)-\wt\xi_{n+1,s+1}^\nu(z))\nonumber\\
 &-b_\l^\nu h_{n+1,s+1}^{\l,\nu}(z+\l)\xi_{n+1,0}^\nu(z+\l)
 +h_{n+1,s+1}^{\l,\nu}(z)\xi_{n+1,0}^\nu(z)\nonumber\\
 &=\left(g_{n+1,s+1}^{\l,\nu}(z)+
 \sum\limits_{i=1}^s B_{i,n-s+i}^\l\xi_{n+1,s+1-i}^\nu(z)\right)
 +(E_{n+1,s+1}^{\l,\,\nu}-f_{n+1,s+1}^{\l,\nu}(z))\xi_{n+1,0}^\nu(z)
 \nonumber\\
 &=E_{n+1,s+1}^{\l,\,\nu}\xi_{n+1,0}^\nu(z)
 +\sum\limits_{i=1}^s B_{i,n-s+i}^\l\xi_{n+1,s+1-i}^\nu(z).
\end{eqnarray}
Now we have constructed a function $\zeta$ having the correct
quasi-periodicity properties (though the first coefficient depends
on $\nu$, so we'll need to deal with this below) but we still cannot
take it to be the function $\xi_{n+1,s+1}^\nu$ that we are trying to
define, as it may not satisfy the equation (\ref{xi.all.s}). We thus
define $R_{n+1,s+1}^{\,\nu}$ to be the ``error'' obtained by
plugging $\zeta$ into (\ref{xi.all.s}):
\begin{equation}\label{R}
  R_{n+1,s+1}^{\,\nu}(z)\xi_{n+1,0}^\nu(z+V):=
  C\,\zeta_{n+1,\, s+1}^{\nu}(z+V)-u_{n}^\nu(z)\,(\zeta_{n+1,\,s+1}^{\nu+1}(z)-
  C\,\xi_{n,\,s}^{\nu+1}(z+V))+\xi_{n,s}^{\nu}(z).
\end{equation}
Notice that for this to make sense we need to assume that we have
been doing all of the above computations simultaneously for $\nu$
and $\nu+1$, so that indeed both $\zeta$'s above are defined at this
point.

From the previous lemma we know that the r.h.s of this formula has
no pole at $\TC_n^\nu$ and vanishes at $\TC_{n}^{\nu+1}-V$, and
thus $R_{n+1,s+1}^{\,\nu}$ is a holomorphic function of $z$. We can
use (\ref{bloch2},\ref{bl4a}) to compute the transformation
properties of $R$ under a shift by a vector $\l\in\Lambda_0$.
Indeed, using (\ref{mon1}) to compute $b_\l^\nu
u_n^\nu(z+\l)=u_n^\nu(z)b_\l^{\nu+1}$, and using (\ref{xi.0.mon})
for the l.h.s., we get, shifting by $\l$ and multiplying by
$b_\l^\nu$, and subtracting the original function,
\begin{eqnarray}\label{R.shift}
 \left(R_{n+1,s+1}^{\,\nu}(z+\l)-R_{n+1,s+1}^{\,\nu}(z)\right)
 \xi_{n+1,0}^\nu(z+V)\\
 =CE_{n+1,s+1}^{\l,\,\nu}\xi_{n+1,0}^\nu(z+V)
 +\sum_{i=1}^s B_{i,n-s+i}^\l\xi_{n+1,s+1-i}^\nu(z+V)\nonumber\\
 -u_n^\nu(z)\left(E_{n+1,s+1}^{\l,\,\nu+1}\xi_{n+1,0}^{\nu+1}(z)
 +\sum_{i=1}^s B_{i,n-s+i}^\l\xi_{n+1,s+1-i}^{\nu+1}(z)\right.\nonumber\\
 \left.-C\sum_{i=1}^s B_{i,n-s+i}^\l\xi_{n,s-i}^{\nu+1}(z+V)\,
 \right)+\sum_{i=1}^s B_{i,n-s+i}^\l\xi_{n,s-i}^{\nu}(z).\nonumber
\end{eqnarray}
Now note that for each constant $B_{i,n-s+i}^\l$ in the above
expression the function it multiplies is exactly the r.h.s.~of
(\ref{xi.all.s}) for $n$ and some $j\le s$, and thus vanishes
identically (this uses in a crucial way the fact that $B$'s do not
depend on $\nu$). Using the formulas (\ref{udd1},\ref{w4}) for
$u_n^\nu$ and $\xi_{n+1,0}$, we get
$$
  R_{n+1,s+1}^{\,\nu}(z+\l)-R_{n+1,s+1}^{\,\nu}(z)=
  C(E_{n+1,s+1}^{\l,\nu}-E_{n+1,s+1}^{\l,\nu+1}).
$$
Moreover, by (\ref{bl5}) we know that the $E'$s are linear functions
of $\l$, i.e. that
$$
  E_{n+1,s+1}^{\l,\nu}-E_{n+1,s+1}^{\l,\nu+1}=2\ell_{n+1,s+1}^\nu(\l)
$$
for some linear function $\ell$; note that
$\ell_{n+1,s+1}^\nu(z)=-\ell_{n+1,s+1}^{\nu+1}(z)$. It then follows
that the difference $R-2\ell$ is periodic with respect to shifts by
$\Lambda_0$, and is thus constant, i.e. we have then
$R_{n+1,s+1}^\nu(z)=2C\ell_{n+1,s+1}^\nu(z)+2A^\nu$. We can now
introduce one last correction and finally define
\begin{equation}\label{sol}
  \xi_{n+1,s+1}^\nu(z):=\zeta_{n+1,s+1}^\nu(z)-(
  \ell_{n+1,s+1}^\nu(z-V/2)+A^\nu+l(z)) \xi_{n+1,0}^\nu(z),
\end{equation}
where $l(z)$ is a linear function such that
$l(V)=A^\nu+A^{\nu+1}$. These functions are going to be solutions
of (\ref{xi.all.s}); indeed, the new error term is equal to
$$
  R_{n+1,s+1}^\nu(z)\xi_{n+1,0}^\nu(z+V)-(\ell_{n+1,s+1}^\nu(z+V/2)
  +A^\nu+l(z+V))\xi_{n+1,0}^\nu(z+V)
$$
$$
  +u_n^\nu(z)(\ell_{n+1,s+1}^{\nu+1}(z-V/2)+A^{\nu+1}+l(z))\xi_{n+1,0}^{\nu+1}(z)
$$
$$
  =\xi_{n+1,0}^\nu(z+V)(R_{n+1,s+1}^\nu(z)-\ell_{n+1,s+1}^\nu(z+V/2)-
  A^\nu-l(z+V)+\ell_{n+1,s+1}^{\nu+1}(z-V/2)+A^{\nu+1}+l(z))
$$
$$
 =\xi_{n+1,0}^\nu(z+V)(2\ell_{n+1,s+1}^\nu(z)+2A^\nu-\ell_{n+1,s+1}^\nu(z)
 -A^\nu-l(V)-\ell_{n+1,s+1}^\nu(z)+A^{\nu+1})=0,
$$
where we used the definitions (\ref{udd1}),(\ref{v0}) and
definitions of $\ell$, $l$, and $A$.

We now need to check that the functions $\xi$ satisfy the
quasi-periodicity conditions (\ref{bloch2}). From
(\ref{xi.0.mon},\ref{bl4a}) it follows that
$$
 b_\l^\nu\xi_{n+1,s+1}^\nu(z+\l)-\xi_{n+1,s+1}^\nu(z)
 =\left(E_{n+1,s+1}^{\l,\nu}-\ell_{n+1,s+1}^\nu(\l)-l(\l)\right)\xi_{n+1,0}^\nu(z)
$$
$$
 +\sum_{i=1}^s B_{i,n-s+i}^\nu\xi_{n+1,s+1-i}^\nu(z),
$$
which means that the function $\xi_{n+1,s+1}^\nu$ satisfies the
quasi-periodicity condition (\ref{bloch2}), if we take
$$
  B_{n+1,s+1}^\l:=E_{n+1,s+1}^{\l,\nu}-\ell_{n+1,s+1}^\nu(\l)-l(\l)
=\frac{E_{n+1,s+1}^{\l,\nu}+E_{n+1,s+1}^{\l,\nu+1}}{2}-l(\l)
$$
(notice that this does not depend on $\nu$, as required in formula
(\ref{bloch2}). Observe that the $B$ we construct is going to depend
on the choice of the linear function $l(\l)$. We have thus
constructed a quasi-periodic solution for $s+1$ and proved the
inductive step of the proposition. \epf

\begin{cor}
For $\xi_{n,s}^\nu$ and $\xi_{n,s}^{\nu+1}$ fixed, the solutions of
(\ref{xi.all.s}), for both $\nu$ and $\nu+1$, are unique up to the
transformation
\begin{equation}\label{ambiguity}
\xi_{n+1,s+1}^\nu(z)\longmapsto\xi_{n+1,s+1}^\nu(z)+ (c+l(z))\xi_{n+1,0}^\nu(z)
\end{equation}
where $c$ is a constant, and $l$ is a linear function on $\C$ such
that $l(V)=0$, both of them independent of $\nu$.
\end{cor}
\bpf This follows by tracing the ambiguity of the choices involved
in the proof of the above lemma. Alternatively, one can prove this
directly by investigating the quasi-periodicity properties of the
difference of two solutions of (\ref{xi.all.s}). \epf

To fix the freedom of choosing $\xi_{n+1,s+1}^\nu$, we would now
like to fix the quasi-periodicity condition to be the same for all
$n$, and to be as simple as possible. Similarly to the case of a
non-degenerate trisecant treated in \cite{kr-tri}, there may be a
problem here in that the functions $\xi_{n,s}^\nu$ may turn out to
be periodic (in our case by ``periodic'' we should mean
$b_{\l_j}\xi_{n,s}^\nu(z+\l_j)=\xi_{n,s}^\nu(z)$). Similarly to the
situation in that paper, note that the space of periodic functions
with a pole on the divisor $\TC_n^\nu$ is the space of sections of
some line bundle, and thus finite-dimensional. Since all divisors
$\TC_n^\nu$ differ by shifts, there is an upper bound on this
dimension independent of $n$ and $\nu$.

It then follows that the functions $\xi_{n,s}^\nu$, for $n$ fixed,
and $s$ and $\nu$ varying, are linearly independent. Indeed, suppose
that there were some linear relation among them, with the maximal
value of $s$ involved in this relation being equal to $S$. But then
solving equations (\ref{xi.all.s}) with $\nu$ and with $\nu+1$,
allows one to express $\xi_{n,S}^\nu$ in terms of
$\xi_{n-1,S-1}^\nu$ and $\xi_{n-1,S-1}^{\nu+1}$, and thus obtain a
linear relation among the $\xi$'s with index $n-1$, and maximal $s$
being equal to $S-1$. By downward induction, we can get to $S=0$ and
get a contradiction with the fact that $\xi_{n,0}^\nu\ne 0$ and is
not proportional to $\xi_{n,0}^{\nu+1}$. Note, moreover, that if for
some $s$ the function $\xi_{n,s}^\nu$ is not periodic, this would
mean that some $B$ is non-zero, and thus $\xi_{n,s+i}^\nu$ could not
be periodic for any $i>0$, as the term in (\ref{bloch2}) with this
non-zero $B$ would be linearly independent with all the other terms
on the right-hand-side there.

\begin{lem}\label{fix.monodromy}
Let $\l_0,\l_1,\ldots,\l_d$ be a set of $\bbC$-linear independent vectors
in $\Lambda_0$. Suppose equations (\ref{xi.all.s}) have periodic
solutions for $i<r$ (and any $n$ and $\nu$), i.e. that there are some
$\Xi_{n,i}^\nu(z)$ such that
$$
 b_{\l_j}^\nu\Xi_{n,i}^\nu(z+\l_j)-\Xi_{n,i}^\nu(z)=0
$$
for all $i<r,$ all $n$ and $\nu$, and such that
$\Xi_{n,0}^\nu=\xi_{n,0}^\nu$ is given by (\ref{w4}). Suppose also
that there are quasi-periodic solutions $\Xi_{n,r}^\nu$ with
\begin{equation}\label{new10}
  b_{\l_j}^\nu\Xi_{n,r}^\nu(z+\l_j)-\Xi_{n,r}^\nu(z)=A_j\xi_{n,0}^\nu(z)
  \qquad \forall j=0,\ldots,d,
\end{equation}
for all $n$, where $A_j$ are some constants such that there does not
exist a linear form $l$ on $\C$ with $l(\l_j)=A_j$, and $l(V)=0$
(i.e. such that the scalar product of the vector
$\vec{A}=(A_1,\ldots, A_d)$ and $V$ is non-zero). Then for all $s\ge
r$, and all $n$ and $\nu$ equations (\ref{xi.all.s}) have
quasi-periodic solutions satisfying (\ref{bloch2}) with
$B_{i,n}^{\l_j}=A_j\delta_{i,r}$, i.e. there exist functions
$\xi_{n,s}^\nu(z)$ for all $s\ge r$, all $n$ and $\nu$ such that
\begin{equation}\label{new20}
  b_{\l_j}^\nu\xi_{n,s}^\nu(z+\l_j)-\xi_{n,s}^\nu(z)=A_j\xi_{n,s-r}^\nu(z).
\end{equation}
(Note that we do not necessarily have
$\xi_{n,i}^\nu(z)=\Xi_{n,i}^\nu(z)$ for $i\le r$, but they satisfy
the same quasi-periodicity, and solve the same equation
(\ref{xi.all.s}).) Moreover, such $\xi_{n,s}^\nu(z)$ are unique up
to adding $c_{n,s}\xi_{n,0}^\nu(z)$, with $c_{n,s}$ being a constant
dependent only on the remainder of $n$ modulo $r$.
\end{lem}
\bpf We prove the lemma by induction in $s$, starting with $s=0$,
with the inductive assumption being that functions $\Xi_{n,i}^\nu$
satisfying (\ref{xi.all.s}) and the quasiperiodicity condition
(\ref{new20}) have been constructed for all $n$ and $\nu$, for all
$i\le r+s$, that they are ``periodic'' for $i<r$, and that moreover
$\xi_{n,i}^\nu(z):=\Xi_{n,i}^\nu(z)$ for $i\le s$ (so that the
inductive assumption for $s=0$ is the assumption of the lemma).

From (\ref{ambiguity}) we know that there must exist solutions
$\wt\xi_{n,s+r+1}^\nu(z)$ of (\ref{xi.all.s}) for all $n$ and $\nu$,
with quasi-periodicity
\begin{equation}\label{extra.term}
  b_{\l_j}^\nu\wt\xi_{n,s+r+1}^\nu(z+\l_j)=
  A_j\Xi_{n,s+1}^\nu(z) +B_{n,s+r+1}^{\l_j}\xi_{n,0}^\nu(z),
\end{equation}
where $B$ are some new constants. The idea now is that we will
adjust all the $\Xi_{n,s+i}^\nu$ for $0<i\le r$ to another set of
solutions of (\ref{xi.all.s}) with the same quasiperiodicity, so
that $\Xi_{n,s+r+1}$ satisfying the quasi-periodicity condition
(\ref{new20}) would exist.

Indeed, suppose we take
$\xi_{n+1,s+1}^\nu(z):=\Xi_{n+1,s+1}^\nu(z)+c_{n+1,s+1}\xi_{n+1,0}^\nu(z)$
for some constant $c_{n+1,s+1}$, independent of $\nu$ (if we added
$l(z)\xi_{n+1,0}^\nu(z)$, the quasi-periodicity of
$\xi_{n+1,s+1}^\nu(z)$ would no longer be the same as that of
$\Xi_{n+1,s+1}^\nu(z)$). If we make such a change, we also need to
add something (let's call it $f^\nu(z)$), to $\Xi_{n+2,s+2}^\nu(z)$,
so that (\ref{xi.all.s}) is still satisfied. Since the $\Xi$'s
themselves satisfied (\ref{xi.all.s}), the corrections we introduce
must also satisfy it, i.e. we must then have
$$
 Cf^\nu(z+V)-u_{n+1}^\nu(z)\left(f^{\nu+1}(z)-Cc_{n+1,s+1}\xi_{n+1,0}^{\nu+1}
 (z+V)\right)+c_{n+1,s+1}\xi_{n+1,0}^\nu(z)=0,
$$
and the same for $\nu+1$. However, this is exactly the equation
(\ref{xi.all.s}) that is satisfied by $c_{n+1}\Xi_{n+2,1}^\nu(z)$,
and thus it follows that $f^\nu(z)=c_{n+1,s+1}\Xi_{n+2,1}^\nu(z)$
would work. Similarly we need to add
$c_{n+1,s+1}\Xi_{n+i+1,i}^\nu(z)$ to each
$\Xi_{n+i+1,s+i+1}^\nu(z)$, so that all of the equations
(\ref{xi.all.s}) are satisfied. Finally in this way we will see that
the necessary adjustment of $\wt\xi_{n+r+1,s+r+1}^\nu$ will be
$$
 \Xi_{n+r+1,s+r+1}^\nu(z):=\wt\xi_{n+r+1,s+r+1}^\nu(z)+
 c_{n+1,s+1}\Xi_{n+r+1,r}^\nu(z)+l_{n+r+1}(z)\xi_{n+r+1,0}^\nu(z),
$$
where we will now need to allow the presence of a linear term to
make the quasi-periodicity be (\ref{new20}) as desired. From
(\ref{extra.term}) and (\ref{new10}) we can compute the
quasi-periodicity to be
$$
 b_{\l_j}^\nu\Xi_{n+r+1,s+r+1}^\nu(z+\l_j)-\Xi_{n+r+1,s+r+1}^\nu(z)
$$
$$
  = A_j\Xi_{n+r+1,s+1}^\nu(z) +\left(B_{n+r+1,s+r+1}^{\l_j}
  +c_{n+1,s+1}A_j+l_{n+r+1}(\l_j)\right)\xi_{n+r+1,0}^\nu(z)
$$
$$
 =A_j\xi_{n+r+1,s+1}^\nu(z) +\left(B_{n+r+1,s+r+1}^{\l_j}
 +(c_{n+1,s+1}-c_{n+r+1,s+1})A_j+l_{n+r+1}(\l_j)\right)\xi_{n+r+1,0}^\nu(z).
$$
For this to be the desired property (\ref{new20}) we must have
$$
 B_{n+r+1,s+r+1}^{\l_j}+(c_{n+1,s+1}-c_{n+r+1,s+1})A_j+l_{n+r+1}(\l_j)=0
 \quad\forall j=0,\ldots, d.
$$
For fixed $n$, this is a system of linear equations for the
difference of the constants $c_{n+1,s+1}-c_{n+r+1,s+1}$ and the
coefficients of the linear form $l$. Recall that $l$ can be chosen
arbitrary such that $l(V)=0$, i.e. if $\l_0\neq 0$ then the
coefficients of $l$ span the $(d)$-dimensional space, in which by
assumption $\vec{A}$ does not lie. Thus the rank of the matrix of
coefficients is $d+1$, and this system of $d+1$ linear equations has
a unique solution. If $\l_0=0$ then the dimension of linear forms
$l$ is $d$, but periodicity condition for $\l_0$ is trivially
satisfied. The inductive assumption is thus proven; note that as a
result we are able to fix the differences
$c_{n+1,s+1}-c_{n+r+1,s+1}$, and thus the constants only depend on
the remainder of $n$ modulo $r$. \epf

\smallskip
\noindent{\bf From local to global considerations.} Until this
point, we have only been working on $\C$, under the assumption that
for all $n$ the functions $\tau_n^\nu(z)$ do not vanish identically.
For $\tau_n^\nu$ given by (\ref{tauexplicit}) that is equivalent to
the assumption that $Un\notin \Sigma$ for all $n$. We now observe
that if a vector $Z\in \bbC^g$ is such that $Z+Un\notin \Sigma$ for
all $n$, then by the same arguments we can construct wave solutions
along the shifted affine subspaces $Z+\C\subset\bbC^d$. Since all
the constructions are explicitly analytic, if we perturb $Z$ (while
still staying away from $\Sigma-Un$), the solutions constructed
along $Z+\C$ will change holomorphically with $Z$. Of course such
solutions can only be constructed locally, while globally there may
be a choice involved, and we may thus have a monodromy for this
choice as we go around $\Sigma-Un$. Thus we cannot a priori expect
$\xi_{n,s}(Z+z)$  (for $z\in\C, Z\in\bbC^g$ to be a global
holomorphic function of $Z$.

Note that for fixed $n$ the functions $\xi_{n+1,s+1}^\nu(Z+z)$ exist
if $Z+nU\notin\Sigma$, and $\xi_{n-i,s-i}^{\nu}(Z+z)$ exist for
$0\le i\le s$. Let us pass now from local to global setting. In this
setting the recurrent equation (\ref{xi.all.s}) takes the form
\begin{equation}\label{xi.all.s1}
  C\,\xi_{s+1}^{\nu}(Z+U+V)-u^\nu(Z)\,(\xi_{s+1}^{\nu+1}(Z+U)-
  C\,\xi_{s}^{\nu+1}(Z+V))+\xi_{s}^{\nu}(Z)=0.
\end{equation}
with
\begin{equation}\label{ff3}
  u^{\nu}(Z)=C{\tau^{\nu+1}(Z+U)\,\tau^{\nu+1}(Z+V)\over
  \tau^{\nu}(Z+U+V)\,\tau^\nu(Z)},
\end{equation}
where
\begin{equation}\label{taunew}
  \tau^\nu(Z)=\theta\left(Z+(1-\nu)W\right)\,
  \left(c_1^{(\,l_1,\,Z)}c_2^{(l_2,Z)}\right)^{\nu-\frac12},
\end{equation}
and $l_1,l_2$ are vectors such that $l_1(V)=l_2(U)=1,\
l_1(U)=l_2(V)=0$. In these notations the arguments in the proof of
proposition \ref{solutions}  yield
\begin{prop}\label{Z.solutions}
If equations (\ref{cm7d}) are satisfied, then:
\begin{itemize}
\item[(i)] for $Z\notin \cup_{i=0}^{N} (\Sigma-iU)$ there exist
functions $\tau_s^{\nu}(Z+z),\, 0\leq s\leq N$, which are local
holomorphic function of $Z$ and  global holomorphic function of
$z\in \C$, such that equations (\ref{xi.all.s1}) hold for
$\xi_s^\nu(Z)=\tau_s^{\nu+1}(Z)/\tau^\nu(Z)$, with
$\tau_0^\nu(Z)=\tau^\nu(Z-U)$ (this is (\ref{w4}).

\item[(ii)] The functions $\xi_s$ satisfy the monodromy relations
\begin{equation}\label{bloch21}
 b_\l^\nu \xi_{s}^\nu (Z+z+\l)-\xi_{s}^{\nu}(Z+z)=
\sum_{i=1}^s B_i^\l(Z)\, \xi_{s-i}^{\nu}(Z+z)\,, \ \l\in \Lambda_0
\end{equation}

\item[(iii)] If $\xi_{s-1}$ is fixed then $\xi_s$ is unique up to the transformation
\begin{equation}\label{transnew}
  \xi_{s}(z+Z)\longmapsto \xi_s(Z+z)+(c_s(Z)+l_s(Z,z))\xi_0,
\end{equation}
where $l_s(Z,z)$ is a linear form in $z$ such that $l_s(Z,V)=0$.
\end{itemize}
\end{prop}

\begin{lem}\label{patch}
Let $r$ be the minimal integer such that
$\xi_1^{\nu},\ldots,\xi_{r-1}^{\nu}$ are ``periodic'' functions of
$z$ with respect to $\Lambda_0$, and such that there is no periodic
solution $\xi_r^{\nu}$ of (\ref{xi.all.s1}). Then the inductive
assumptions of lemma (\ref{fix.monodromy}) are satisfied, i.e. the
quasi-periodicity coefficients $B_i^\l(Z)$ in (\ref{bloch21}) do not
depend on $Z$ or $i-s$.
\end{lem}
\bpf
By assumption $\xi_{r-1}^\nu(z)$ is ``periodic'', i.e. we have
$$
 0= b_\l^\nu \xi_{r-1}^\nu (Z+z+\l)-\xi_{r-1}^{\nu}(Z+z)=
\sum_{i=1}^{r-1} B_i^\l(Z)\, \xi_{r-1-i}^{\nu}(Z+z)
$$
However, as it was noted above, the functions $\tau_s^{\nu}$, for
$0\le s\le r-1$ (recall that if not, by applying (\ref{xi.all.s1})
we could produce a linear dependence having only one term, which is
impossible), which means that all the coefficients $B_i^\l$ are zero
for all $i\le r-1$. Thus the monodromy of the next function is given
by
$$
  b_\l^\nu\xi_r^\nu(Z+z+\l)-\xi_r^\nu(Z+z)
 =B_r^\l(Z)\xi_0^\nu(Z+z),
$$
where we of course know $\xi_0^\nu$ explicitly, and $B_r^\l$ is a
local function of $Z$ defined locally for
$$
 Z\in X\setminus\bigcup_{i=0}^{r-1}(\Sigma-iV)
$$
From lemma \ref{Z.solutions} we know that the only ambiguity in the
choice of the solutions $\xi_{n,r}^\nu(z)$ is given by
(\ref{transnew}). Recall that adding a linear function multiple will
change the equation to be satisfied, while adding a constant
multiple does not change the quasi-periodicity properties, so that
finally $B_r^\l(Z)$ is independent of the ambiguity, and is
well-defined as a holomorphic function of $Z\in X'$. Since the locus
$\Sigma\subset X$ is of codimension at least 2, by Hartogs theorem
the function $B_r^\l(Z)$ can be extended holomorphically to all of
$X$. Since $X$ is compact, this means that $B_r^\l(Z)$ is a
constant, which we can denote $A_\l$ for the inductive assumption of
lemma \ref{fix.monodromy}. If we had $\vec A\cdot V=0$, then by a
transformation (\ref{ambiguity}) with a suitable linear term we
could get a new solution with $A_{\l_i}=0$ for $i=0...d$, i.e. the
function $\xi_{n,r}^\nu(z)$ could be made ``periodic'',
contradicting the way we chose $r$. \epf

\begin{lem}\label{no.sigma}
In the setup of our construction the ``bad locus'' $\Sigma$ is
actually empty, i.e. if equation (\ref{cm7d}) (part (C), the weakest
assumption of our main theorem) is satisfied, then
$\Sigma=\emptyset$.
\end{lem}
\bpf The proof of this lemma is analogous to the proof of the
similar statement for the fully discrete trisecant characterization
of Jacobians treated in \cite{kr-tri}, once we first prove that
$\Sigma^0=\Sigma^1$.

The only ambiguity in the definition of $\tau_1^\nu(Z)$ is in the
choice of the coefficient $c_1$ in (\ref{transnew}). Suppose there
exists a point $A\in\Sigma^0\setminus\Sigma^1$, i.e. such that
$\theta(A+NV)\ne0=\theta(A+NV+W)$ for some $N$. Then locally near
the point $A+NV$ choose some holomorphic branch of the function
$\xi_1^1(Z)=\tau_1^0(Z)/\tau^1(Z)$. Doing this fixes the value of
$c_1(Z)$ for all $Z$ near $A+NV$. However, since the ambiguity in
the choice of $\xi_1^0(Z)=\tau_1^1(Z)/\tau^0(Z)$ is given by the
same function $c_1(Z)$ (which did not depend on $\nu$!) it means
that for $Z$ in a neighborhood of $A+NV$, but away from $\Sigma^0$,
we also have a fixed choice of $\xi_1^0(Z)$, and thus also of the
holomorphic function $\tau_1^1(Z)$. Since $\Sigma^0$ has codimension
at least 2 in $X$, the function $\tau_1^1(Z)$ can thus be extended
to all points in a neighborhood of $A$, which is a contradiction.
Thus we must have $\Sigma^0\subset\Sigma^1$, and of course by
symmetry they are in fact equal.

We now prove in the same manner that $\Sigma=\Sigma+rU$: above we
used the fact that in (\ref{transnew}) $c_1$ is independent of
$\nu$, and now we will use $c_1(Z)=c_1(Z+rU)$. Indeed, suppose we
have $A\in(\Sigma-rU)\setminus\Sigma$. This means that in a
neighborhood of $p$ we can choose locally holomorphically the
function $\tau_1^\nu(Z)$, i.e. chose a local holomorphic branch of
$c_1(Z)$. However, since $c_1(Z)=c_1(Z+rU)$, this also fixes the
choice of $c_1$ in a neighborhood of the point $p+rU$, and thus in a
neighborhood of $p$, outside of $\Sigma-rU$, we have a holomorphic
function $\tau_1^\nu(Z+rU)$, which now can be extended across
$\Sigma-rU$, which we know to be of codimension at least two. This
constructs a solution $\tau_1^\nu(A+rU)$, which contradicts the
assumption $A\in(\Sigma-rU)$.

Since by definition $\Sigma$ has no subset invariant under a shift
by $U$, this implies that either $\Sigma$ is empty, or $r>1$.
Suppose now that $\Sigma$ is non-empty, so $r>1$. Recall that we
have $\tau_0^\nu(Z)=\tau^\nu(Z-U)=\theta(Z+(1-\nu)W-U){\rm const}$,
and thus since $\Sigma^0=\Sigma^1$ we have
$\tau_0^\nu|_{\Sigma+U}=0$. Thus for any $Z\in\Sigma+U$, for $s=1$
the last two terms in (\ref{xi.all.s1}) vanish, yielding
$$
 C\xi_1^\nu(Z+U+V)-u^{\nu}(Z)\xi_1^{\nu+1}(Z+U)=0
$$
However, this is exactly equation (\ref{xi.all.s1}) for $s=0$, which
is solved by $\xi_0^\nu$, and thus all {\it periodic with respect to
$\Lambda_0$} solutions are constant multiples of $\xi_0^\nu$. By
using (\ref{transnew}) we can subtract this constant, and get a
solution such that $\xi_1^\nu(Z+U)=0$ for any $Z\in\Sigma+U$, i.e.
we have $\xi_1^\nu|_{\Sigma+2U}=0$.

Now we can repeat this process: indeed, for $s=2$ and
$Z\in\Sigma+2U$ the last two terms in (\ref{xi.all.s1}) have
$\xi_1^\nu(Z+V)$ and $\xi_1^\nu(Z)$ appearing, and thus vanish, so
that as a result we see that $\xi_2^\nu$ on $\Sigma+3U$ is a
constant multiple of $\xi_0^\nu$. By using (\ref{transnew}) again,
we can make this multiple to be zero again. Repeating this a number
of times, we will eventually get $\xi_{r-1}^\nu|_{\Sigma+rU}=0$.

Since $\Sigma=\Sigma+rU$, we also have $\tau^\nu|_{\Sigma+rU}=0$,
and thus for $Z\in\Sigma+rU$ and $s=r-1$ the last two terms in
(\ref{xi.all.s1}) vanish to the second order --- both factors of
each summand vanish. Thus $\xi_r^\nu$ can be defined in a
neighborhood of $\Sigma+rU=\Sigma$ as a holomorphic function
vanishing on $\Sigma+rU$. However, this implies in particular that
$b_\l^\nu\xi_r^\nu(Z+\l_j)-\xi_r^\nu(Z)=0$ for $Z\in\Sigma$ and any
$\l_j\in\Lambda_0$, which contradicts the assumption that
$\xi_r^\nu$ could not be periodic. The lemma is thus proven.
\epf

As shown above, if $\Sigma$ is empty, then the functions
$\tau_s^{\nu}$ can be defined as global holomorphic functions of
$Z\in \bbC^g$. Then, as a corollary of the previous lemmas we get
the following statement.

\begin{lem} Suppose (\ref{cm7d}) for $\theta(Z)$ holds.
Then there exists a pair of formal solutions
\begin{equation}\label{gff1}
  \phi^\nu=\sum_{s=0}^{\infty}\xi_s^{\nu}(Z)\,k^{-s}
\end{equation}
of the equation
\begin{equation}\label{gff2}
kC\phi^{\nu}(Z+U+V,k)-u^{\nu}(Z)(k\phi^{\nu+1}(Z+U, k)-
C\phi^{\nu+1}(Z+V,k))-\phi^{\nu}(Z, k)=0\,,
\end{equation}
with $C=c_3$ and
\begin{equation}\label{gff3}
  u^{\nu}(Z)={\tau^{\nu+1}(Z+U)\,\tau^{\nu+1}(Z+V)\over
  \tau^{\nu}(Z+U+V)\,\tau^\nu(Z)},
\end{equation}
where $\tau^{\nu}$ is given by (\ref{taunew}), such that:
\begin{itemize}
\item[(i)] the coefficients $\xi_s^\nu$ of the formal series
$\phi^{\nu}$ are of the form
$\xi_s^{\nu}(Z)=\tau_s^{\nu+1}(Z)/\tau^{\nu}(Z)$, where
$\tau_s^{\nu} (Z)$ are holomorphic functions;
\item[(ii)] $\phi^{\nu}(Z,k)$ is quasi-periodic with respect to the
lattice $\Lambda$ and for the basis vectors $\l_j$ in $\C$ its
monodromy relations have the form
\begin{equation}\label{gff4}
  \phi^{\nu}(Z+\lambda_j)=(1+A_{\l_j}\,k^{-1})\,\phi^{\nu}(Z,k),\ \ j=1,\ldots, g,
\end{equation}
where $A_{\l_j}$ are constants such that there is no linear form
on $\C$ vanishing at $V$, i.e. $l(V)=0$, and such that
$l(\l_j)=A_{\l_j}$;
\item[(iii)] $\phi^{\nu}$ is unique up to the multiplication by a
constant in $Z$ factor.
\end{itemize}
\end{lem}

\section{The spectral curve}
In this section we finish a proof of the fact that condition $(C)$
of the main theorem characterizes Prym varieties. Indeed, in the
previous section we showed that if $(C)$ holds, some quasi-periodic wave
solutions can be constructed. In this section we show that these
wave solutions are eigenfunctions of commuting difference
operators, and identify $X$ with the Prym variety of the spectral
curve of these operators. Much of the argument is analogous to
that in \cite{kr-prym}.

The formal series $\phi^{\nu}(Z,k)$ constructed in the previous
section define a wave function
$$\psi=\psi_{nm}(k):=k^{n}\phi^{\nu_{nm}}(nU+mV+Z,k).$$
This wave function determines a unique pseudo-difference operator
$\L$ such that $\L\psi=k\psi$ (the coefficients of this $\L$ can
be computed inductively term by term); we note that the ambiguity
in the definition of $\phi^{\nu}(Z)$ (it is only defined up to a
factor that is $T_1$-invariant) does not affect the coefficients
of the wave operator. Therefore, its coefficients are of the form
\begin{equation}\label{kkd}
\L=\sum_{s=-1}^{\infty} w_s^{\nu_{nm}}(Z+nU+mV)\,T_1^{-s},
\end{equation}
where $w_s^{\nu}(Z)$ are well-defined meromorphic sections of line
bundles on $X$ with automorphy properties given by (\ref{blambda}).

As before, we define functions $\wt F_j$ by formula (\ref{gen1}), i.e. we set
$$
\wt F_j:= \res_T\left((\L^jT_1^{-1}-T_1\L^j)(T_1-T_1^{-1})^{-1}\right)
$$
The definition of $\psi$ implies that these functions are of the form
\begin{equation}\label{Fj.f}
\wt F_j=\wt F_j^{\,\nu_{nm}}(Un+Vm+Z)
\end{equation}
where $\wt F_j^{\,\nu}(Z)$ are meromorphic functions on $X$.

\begin{lem} There exist vectors $V_{m}=\{V_{m,k}\}\in\bbC^g$
and constants $v_{m}\in\bbC$ such that
\begin{equation}\label{nnov7}
  \wt F_{j}^\nu (Z)=v_{j}+ \frac{\p}{\p V_j}\left(\ln \tau^{\nu}(Z)-\ln \tau^{\nu+1}(Z+U)\right).
\end{equation}
\end{lem}
\bpf
Consider the formal series $\psi^\s$ given by (\ref{more5}). It has the form
\begin{equation} \label{501}
  \psi^\s_{n,m}=k^{-n}\phi^{\,\s,\,\nu_{n,m}}(Un+Vm+Z,k),
\end{equation}
where the coefficients of the formal series
\begin{equation}\label{502}
  \phi^{\,\s,\,\nu}(Z,k)=\sum_{s=0}^{\infty}\xi_s^{\,\s,\,\nu}(Z)\,k^{-s}
\end{equation}
are difference polynomials in the coefficients of $\phi^{\nu}$ and
$\phi^{\nu+1}$. Therefore, we know a priori that
$\xi^{\,\s,\,\nu}(Z)$ are meromorphic functions, which may have
poles for $Z\in\TC^\nu, Z\in\TC^{\nu+1}$, or for $Z$ on the
translates of these two divisors by integer multiples of $U$. We
claim that in fact these coefficients are of the form
\begin{equation}\label{503}
\xi_s^{\,\s,\,\nu}={\tau_s^{\,\s,\,\nu}(Z)\over \tau^\nu(Z)},
\end{equation}
where $\tau_s^{\,\s,\,\nu}(Z)$ are some holomorphic functions, i.e.
that they only have simple poles at $\TC^\nu$.

Indeed, we showed in section 3 that $\psi^{\s}$ solves the
equation $H\psi^{\s}=0$. In section 4 we deduced from the
statement $(C)$ of the main theorem the fact that $\psi^\nu$ may
only have a simple pole on $\TC^\nu$. By replacing $\psi^\nu$ by
$\psi^{\,\s,\,\nu}$, and replacing $U$ by $-U$ we get from
statement $(C)$ functional equations for $\tau_s^{\,\s,\,\nu}$ and
in the same way deduce also that $\psi^{\s,\nu}$ only has pole at
$\TC^\nu$.

Equation (\ref{J0}) then implies that $\wt \F_j^{\nu}$ are the
coefficients of the formal series
\begin{equation}\label{J05}
  -k+(k^2-1)\sum_{j=1}^{\infty}\wt F_j^{\,\nu}(Z)=
  k^{-1}\phi^{\,\s,\,\nu+1}(Z+U,k)\,\phi^{\nu}(Z,k)-
  k\phi^{\,\s,\,\nu}(Z,k)\,\phi^{\nu+1}(Z+U,k)
\end{equation}
It thus follows that $\wt F^{\,\nu}_j(Z)$ have simple poles only
at the divisors $\TC^{\nu}$ and $\TC^{\nu+1}-U$ --- these are the
only possible poles of the right-hand-side. Moreover, equation
(\ref{q81}) says (recall that ${\bf t_1}$ is shifting the variable
$n$, i.e. adding $U$) that there exist meromorphic functions
$Q^{\,\nu}_j$ such that
\begin{equation}\label{504}
  \wt F_j^{\,\nu}(Z)=Q_j^{\,\nu}(Z)-Q_j^{\,\nu+1}(Z+U).
\end{equation}
We know a priori that $Q_j^{\,\nu}$ may only have poles at
$\TC^{\,\nu}$ and $\TC^{\,\nu+1}-U$. However, if there were a pole
at $\TC^{\nu+1}-U$, it would then mean that $Q_j^{\nu+1}(Z+U)$
would have a pole at $\TC^\nu-2U$, and since by our initial
assumptions $U$ was not a point of order two, this is impossible.
Thus $Q^{\,\nu}_j$ has simple pole only on $\TC^\nu$, as desired
for the expression (\ref{nnov7}) for $\wt F_j^\nu$ to be valid.
The functions $\wt F_j^{\nu}$ are abelian functions. Therefore,
the residue of $Q_j^\nu$ is a well-defined section of the
theta-bundle restricted $\TC^{\,\nu}$, i.e.
$$\left(Q_j^{\,\nu}\tau^{\nu}\right)|_{\TC^{\,\nu}}\in
H^0(\tau^{\nu}|_{\TC^{\,\nu}})$$
It is know that the later space is spanned by
the directional derivatives of the theta function. Thus we see that
there must exist some vector $V_j^\nu\in \bbC^g$
such that $Q_j^\nu-\left(\p\ln\tau^{\nu}(Z)/\p{V_j^{\,\nu}}\right)$ is a holomorphic
function. The periodicity of $\wt F_j^\nu$ with respect to the lattice
implies that $V_j^\nu=V_j^{\nu+1}$, and thus (\ref{nnov7}) holds.
\epf

Consider now the linear space spanned the functions $\{\wt
F_j^{\nu}(Z), \, j= 1,\ldots\}$. From (\ref{nnov7}) we see that
there are only $g+1$ parameters involved in determining $\wt
F_j^{\nu}$, and thus this space is at most $g+1$-dimensional.
Therefore, for all but $\wt g:=\dim\ \{\wt F_j^{\nu}(Z)\}-1\leq g$
positive integers $j$, there exist constants $c_{i,j}$ such that
\begin{equation}\label{f1}
  \wt F_j^{\nu}(Z)=c_{0,j}+\sum_{i=1}^{j-1} c_{i,j}\wt F_i^\nu(Z).
\end{equation}
Let $I$ denote the subset of integers $j$ for which there are no
such constants. We call this subset the gap sequence --- the
corresponding $\wt F_j^{\nu}$ form a basis.

\begin{lem}
Let $\L$ be the pseudo-difference operator corresponding to the
quasi-periodic (Bloch) wave function $\psi$ constructed above. Then, for the
difference operators
\begin{equation}\label{a2}
  \wh L_j:=L_j+\sum_{i=1}^{j-1} c_{i,j}L_{n-i}=0, \quad \forall j\notin I,
\end{equation}
the following equations are satisfied with some constants $a_{s,j}$:
\begin{equation}\label{lp}
  \wh L_j\,\psi=a_j(k)\,\psi, \ \ \ a_j(k)=k^j+\sum_{s=1}^{\infty}a_{s,j}k^{j-s}
\end{equation}
\end{lem}
\bpf
From the proof of Theorem 3.5 we get
$$
  [L_j,H]\equiv \left({\bf t}_2 \wt F_j-\wt F_j\right) \, (T_1-T_2)\mod \O_H
$$
Therefore, operators $\wh L_j$ and $H$ commute in $\O/\O_H$.
Hence, if $\psi$ is a Bloch wave solution of
(\ref{laxdd1}), i.e. $H\psi=0$, then $\wh L_j\psi$ is also a
Bloch solution of the same equation. Since (\ref{laxdd1})
has a unique solution up to multiplication by constant (i.e. the
kernel of $H$ is one-dimensional), we must have $\wh
L_j\psi=a_j(Z,k)\psi$, where $a_j$ is $T_1$-invariant, i.e.
$a_j(Z,k)=a_j(Z+U,k)$.

Note that the constant factor ambiguity in the definition of $\psi$
does not affect $a_j$, and thus $a_j$ are well-defined {\it global}
meromorphic functions on $\bbC^g\setminus \Sigma$. Since the closure
of $\bbZ U$ is dense in $X$, the $T_1$ invariance of
$a_j$ implies that $a_j$ is a holomorphic function of
$Z\in X$, and thus it is constant in $Z$ (note that we
in fact have $a_{s,n}=-c_{s,n}$ for $s\leq n$). \epf

If we set now $m=0$, the operator $\wh L_j$ can be regarded as
a ${Z}$-parametric family   of ordinary difference operators $\wh
L_j^{Z}$.
\begin{cor}
The operators $\wh L_j^{Z}$ commute with each other,
\begin{equation}\label{com1}
  [\wh L_i^{Z},\wh L_j^{Z}]=0.
\end{equation}
\end{cor}

A theory of commuting difference operators containing a pair of
operators of co-prime orders was developed in  \cite{mum,kr-dif}.
It is analogous to the theory of rank 1 commuting differential
operators \cite{ch1,ch2,kr1,kr2,mum} (relatively recently this
theory was generalized to the case of commuting difference operators
of arbitrary rank in \cite{n-kr}.)

\begin{lem}
Let $\A^Z$ be the commutative ring of ordinary difference operators
spanned by the operators $\wh L_j^Z$. Then there exists an irreducible
algebraic curve $\G$ of arithmetic genus $\hat g$ with involution
$\s:\G\longmapsto \G$ such that for a generic $Z$ the ring $\A^Z$ is
isomorphic to the ring of meromorphic functions on $\G$ with the
only poles at two smooth points $P_1^{\pm}$, which are odd with
respect to the involution $\s$. The correspondence $Z\to \A^Z$
defines a holomorphic map of $X$ to the space of odd torsion-free
rank one sheaves $\F$ on $\G$
\begin{equation}\label{is}
j: X\longmapsto \overline{\rm Prym}(\G)=Ker (1+\s)\subset \overline{\rm Pic}(\G).
\end{equation}
\end{lem}
{\it Proof.} As shown in \cite{mum,kr-dif}
there is a natural correspondence
\begin{equation}\label{corr}
\A\longleftrightarrow \{\G,P_{\pm},  \F\}
\end{equation}
between commutative rings $\A$ of ordinary linear
difference operators containing a pair of monic operators of co-prime orders, and
sets of algebro-geometric data $\{\G,P_{\pm}, [k^{-1}]_{\pm}, \F\}$,
where $\G$ is an algebraic curve with fixed first jets $[k^{-1}]_{\pm}$ of local
coordinates $k^{-1}_{\pm}$ in the neighborhoods of smooth
points $P_1^{\pm}\in\G$, and $\F$ is a torsion-free rank 1 sheaf on $\G$ such that
\begin{equation}\label{sheaf}
  h^0(\G,\F)=h^1(\G,\F(nP_+-nP_-)=0.
\end{equation}
The correspondence becomes one-to-one if the rings $\A$ are
considered modulo conjugation $\A'=g\A g^{-1}$.

The construction of the correspondence (\ref{corr}) depends on a
choice of initial point $n_0=0$. The spectral curve and the sheaf
$\F$ are defined by the evaluations of the coefficients of
generators of $\A$ at a finite number of points of the form $n_0+n$.
In fact, the spectral curve is independent on the choice of $x_0$,
but the sheaf does depend on it, i.e. $\F$ depends on the choice of $n_0$.

Using the shift of the initial point it is easy to show that the
correspondence (\ref{corr}) extends to the commutative rings of
operators whose coefficients are {\it meromorphic} functions of $x$.
The rings of operators having poles at $n=0$ correspond to sheaves
for which the condition (\ref{sheaf}) for $n=0$ is violated.

A commutative ring $\A$ of linear ordinary difference operators is
called maximal if it is not contained in any larger commutative
ring. The algebraic curve $\G$ corresponding to a maximal ring is
called the spectral curve of $\A$. The ring $\A$ is isomorphic to
the ring $A(\G,P_1^{\pm})$ of meromorphic functions on $\G$ with the
only pole at $P_1^+$, and vanishing at $P_1^-$. The isomorphism is
given by the equation
\begin{equation}\label{z2}
  L_a\psi=a\psi, \ \ L_a\in \A, \ a\in A(\G,P_1^{\pm}),
\end{equation}
where $\psi$ is a common eigenfunction of the commuting operators.

Let $\G^Z$ be the spectral curve corresponding to the maximal ring
$\hat \A^Z$ containing $\A^Z$. The eigenvalues $a_j(k)$ of the
operators $\hat L_j^Z$ defined in (\ref{lp}) coincide with the
Laurent expansions at $P_1^+$ of the meromorphic functions $a_j\in
A(\G^Z,P_{\pm})$, and thus are $Z$-independent. Hence, the spectral
curve $\G^Z$ in fact does not depend on $Z$.

The functions $\psi^\s$ are eigenfunctions of $\wh L_j$:
\begin{equation}\label{505}
  \wh L_j\psi^\s=-a_j(k)\psi^{\s}.
\end{equation}
Hence, the correspondence $\psi\leftrightarrow \psi^\s$ gives rise to an involution $\s$ of the spectral curve. The eigenvalues $a_j$ are
odd with respect to the involution, and thus the lemma is proved. \epf

The next step is to consider deformations of $\A^Z$ defined by the
discrete Novikov-Veselov hierarchy introduced in Section 3. Through this
hierarchy we identified the space spanned by functions $\wt F_j$
with the tangent space to the orbit of the hierarchy. Lemma 5.1
identifies the orbit of the hierarchy with $Z+Y$, where Y is the
closure of the group spanned by vectors $V_j$. The orbit of the NV
hierarchy is the odd part of the orbit of two Kadomtsev-Petviashvili flows
corresponding to points $P_1^{\pm}$. It follows from
\cite{shiota} that the orbit of the discrete NV hierarchy is
isomorphic to $\P(\G)$. For a generic $Z$ the ring $A^Z$ is
a maximal odd ring. Therefore, we get

\begin{lem} For $Z\in\bbC^g$ generic, the orbit of $\A^Z$ under
the NV flows defines an isomorphism:
\begin{equation}\label{imb}
  i_Z:\P(\G)\longmapsto Z+Y\subset X.
\end{equation}
\end{lem}

\begin{cor}
The Prym variety $\P(\G)$ of the spectral curve $\G$ is compact.
\end{cor}
The compactness of the Prym variety is not as restrictive as the
compactness of the Jacobian (see \cite{mdl}). Nevertheless, it
implies an explicit description of the singular points of the
spectral curve.  The following result is due to Robert Friedman (see
the appendix of \cite{kr-prym}):
\begin{cor} The spectral curve $\G$ is smooth outside of fixed points $q_k$ of
the involution $\s$. The branches of $\G$ at $q_k$ are linear and
are not permuted by $\s$.
\end{cor}
The arguments identical to that used at the end of \cite{kr-prym}
prove that in fact the singular points $q_k$ are at most double
points. For a curve with at most double singular points all the
sheafs $\F$ are line bundles. Therefore, the map $j$ in (\ref{is})
is inverse to $i_Z$ in (\ref{imb}), and the main theorem is thus
proven.

\medskip

{\footnotesize{\sl Acknowledgments. \rm We would like to thank Enrico Arbarello,
Robert Friedman and Takahiro Shiota, in very useful discussions with
whom we have learned many things about the geometry of Prym varieties and Prym theta
divisors.}}


\begin{thebibliography}{**}
\bibitem{flex}E. Arbarello, I. Krichever, G. Marini,
{\it Characterizing Jacobians via flexes of the Kummer Variety},
Math. Res. Lett. {\bf 13} (2006) 1, 109--123.

\bibitem{bd}
A. Beauville, O. Debarre, {\it Sur le probl\`eme de Schottky pour
les vari\'et\'es de Prym.} Ann. Scuola Norm. Sup. Pisa Cl. Sci. 
{\bf 14}  (1987) 4, 613--623 (1988).

\bibitem{ch1}
J.L.Burchnall, T.W. Chaundy, {\it Commutative ordinary differential
operators.I}, Proc. London Math. Soc. {\bf 21} (1922), 420--440.

\bibitem{ch2}
J.L.Burchnall, T.W. Chaundy, {\it Commutative ordinary differential
operators. II}, Proc. Royal Soc. London  {\bf 118} (1928), 557--583.

\bibitem{deb}
O. Debarre: {\it Vers une stratification de l'espace des modules des
vari\'et\'es ab\'eliennes principalement polaris\'ees.} Complex
algebraic varieties (Bayreuth, 1990), 71-–86, Lecture Notes in
Math., 1507, Springer, Berlin, 1992.

\bibitem{mdl}
P. Deligne, D. Mumford, {\it The irreducibility of the space of
curves of given genus.} Inst. Hautes Etudes Sci. Publ. Math. {\bf 36} (1969), 75--109.

\bibitem{grin}
A. Doliwa, P. Grinevich, M. Nieszporski, P. M. Santini, {\it
Integrable lattices and their sub-lattices: from the discrete Moutard (discrete Cauchy-Riemann) 4-point equation to the self-adjoint 5-point scheme}, arXiv:nlin/0410046.

\bibitem{dkn}
B.A. Dubrovin, I.M. Krichever, S.P. Novikov {\it Schr\"odinger
equation in magnetic field and Riemann surfaces}, Dokl. Akad. Nauk SSSR {\bf 229} (1976), 1, 15--18.

\bibitem{fay}
J.D. Fay, Theta functions on Riemann surfaces. Lecture Notes
in Mathematics, Vol. 352. Springer-Verlag, Berlin-New York, 1973.

\bibitem{fay2}
J.D. Fay, {\it On the even-order vanishing of Jacobian theta
functions.}  Duke Math. J.  {\bf 51} (1984) 1, 109--132.

\bibitem{gun1}
R. Gunning, {\it Some curves in abelian varieties,} Invent. Math. {\bf 66} (1982) 3, 377--389.

\bibitem{kr1}
I.M. Krichever, {\it Integration of non-linear equations by methods
of algebraic geometry}, Funct. Anal. Appl. {\bf 11} (1977) 1, 12--26.

\bibitem{kr2}
I.M. Krichever, {\it Methods of algebraic geometry in the theory of
non-linear equations}, Russian Math. Surveys {\bf 32} (1977) 6,
185--213.

\bibitem{kr-dif}
I.Krichever, {\it Algebraic curves and non-linear difference equation}, Uspekhi Mat. Nauk {\bf 33} (1978) 4, 215--216.

\bibitem{kr-schot}
I.Krichever, {\it Integrable linear equations and the
Riemann-Schottky problem},  Algebraic geometry and number theory,
497--514, Progr. Math., 253, Birkh\"auser Boston, Boston, MA, 2006.

\bibitem{kr-prym}
I.Krichever, {\it A characterization of Prym varieties.}
Int. Math. Res. Not.  2006, Art. ID 81476, 36 pp.

\bibitem{kr-tri}
I. Krichever, {\it Characterizing Jacobians via trisecants of the
Kummer Variety}, math.AG/0605625.

\bibitem{kr5}
I. Krichever, {\it Two-dimensional periodic difference operators
and algebraic geometry.} (Russian) Dokl. Akad. Nauk SSSR {\bf 285}
(1985) 1, 31--36.

\bibitem{n-kr}
I. Krichever, S. Novikov,
{\it Two-dimensional Toda lattice, commuting difference operators
and holomorphic vector bundles}, Uspekhi Mat. Nauk , {\bf 58} (2003) 3, 51--88.

\bibitem{kp1}
I.Krichever, D.H.Phong, {\it On the integrable geometry of soliton equations and $N=2$ supersymmetric gauge theories.}
J. Differential Geom. {\bf 45} (1997) 2, 349--389.

\bibitem{kwz}
I. Krichever, P. Wiegman, A. Zabrodin, {\it Elliptic solutions to
difference non-linear equations and related many-body problems.}
Comm. Math. Phys.  {\bf 193}  (1998) 2, 373--396.

\bibitem{mum2}
D. Mumford {\it Theta characteristics of an algebraic curve.}
Ann. Sci. \'Ecole Norm. Sup. {\bf 4} (1971), 181--192.

\bibitem{mum1}
D. Mumford, {\it Prym varieties. I.} Contributions to analysis (a
collection of papers dedicated to Lipman Bers),  325--350.
Academic Press, New York, 1974.

\bibitem{mum}
D. Mumford, {\it An algebro-geometric construction of commuting
operators and of solutions to the Toda lattice  equation,
Korteweg-de Vries equation and related non-linear equations}.
Proceedings Int.Symp. Algebraic Geometry, Kyoto, 1977, 115--153,
Kinokuniya Book Store, Kyoto, 1978.

\bibitem{nv}
S. Novikov, A. Veselov, {\it Finite-gap two-dimensional periodic
Schr\"odinger operators: exact formulae and evolution equations},
Dokl. Akad. Nauk SSSR {\bf 279} (1984) 1, 20--24.

\bibitem{wilson}
G.Segal, G.Wilson, {\it Loop groups and equations of KdV type},
Inst. Hautes Etudes Sci. Publ. Math. {\bf 61} (1985), 5--65.

\bibitem{serr}
J-P. Serre, {\it Faisceaux alg\'ebriques coh\'erents},  Ann. of
Math. {\bf 61} (1955), 197--278.

\bibitem{shiota}
T. Shiota, {\it Characterization of Jacobian varieties in terms of
soliton equations}, Invent. Math. {\bf 83} (1986) 2, 333--382.

\bibitem{wel1}
G.E. Welters, {\it A criterion for Jacobi varieties}, Ann. of
Math., {\bf 120} (1984) 3, 497--504.

\end{thebibliography}
\end{document}